\documentclass[12pt,reqno]{amsart}

\usepackage{undertilde, slashbox}
\usepackage[final
]{showkeys}

\usepackage{AKstyle}

\usepackage[table]{xcolor}

\newcommand{\vertiii}[1]{{\left\vert\kern-0.25ex\left\vert\kern-0.25ex\left\vert #1
    \right\vert\kern-0.25ex\right\vert\kern-0.25ex\right\vert}}

\numberwithin{equation}{section}

\newcommand{\fg}{\frak{g}}

\usepackage{caption}
\usepackage[labelfont=rm]{subcaption}

\usepackage[pagebackref=true, colorlinks]{hyperref}

\hypersetup{pdffitwindow=true,linkcolor=blue,citecolor=blue,urlcolor=blue,menucolor=blue}

\usepackage{comment}

\begin{document}

\author{Dubi Kelmer}
\thanks{Kelmer is partially supported by NSF grant DMS-1237412 and NSF grant DMS-1401747.}
\email{kelmer@bc.edu}
\address{Boston College, Boston, MA}
\author{Alex Kontorovich}
\thanks{Kontorovich is partially supported by
an NSF CAREER grant DMS-1254788/DMS-1455705, an NSF FRG grant DMS-1463940, and an Alfred P. Sloan Research Fellowship.}
\email{alex.kontorovich@rutgers.edu}
\address{Rutgers University, New Brunswick, NJ}

\title
{Effective equidistribution of  shears and applications}

\begin{abstract}
A ``shear'' is a unipotent translate of a cuspidal geodesic ray in the quotient of the hyperbolic plane by a non-uniform discrete subgroup of $\PSL(2,\R)$, possibly of infinite co-volume. We prove the regularized equidistribution of shears under large translates with effective (that is, power saving) rates. We also give applications 
to weighted second moments of GL(2) automorphic L-functions, and to counting lattice points on locally affine
symmetric spaces.
\end{abstract}
\date{\today}
\maketitle
\tableofcontents

\section{Introduction}\label{sec:intro}

In this paper, we
 prove the effective  (meaning, with power savings rate)  equidistribution of ``shears'' (see below for definitions) of ``cuspidal'' geodesic rays on hyperbolic surfaces.
 Our proofs are
  quite ``soft,'' in that
  we only use
 mixing
and 
standard
properties of Eisenstein series, 
 rather than
explicit
 spectral decompositions, special functions, or any estimates on time spent near a cusp.
 This allows us to 
 extend our methods to surfaces of infinite volume (in fact the proofs are
 surprisingly
  easier in this case).
 %
As a direct consequence,
we
complete the 
general
problem
of 
obtaining 
 effective
 asymptotics for
 counting (in archimedean norm balls)
discrete orbits on affine quadrics; as discribed in \secref{sec:taxonomy}, exactly two lacunary settings remained unsolved, which are settled in this paper.
Another application is to weighted second moments of $\GL(2)$ automorphic $L$-functions.

When the 
surface
has
infinite volume,
we discover 
two new and 
completely unexpected
 phenomena: (1) the orbit count asymptotic
can be proved with a
{\it uniform}
power savings error in congruence towers {\it without} inputting any information on the spectral gap.\footnote{By ``spectral gap'' we always mean the distance between the first eigenvalue $\gl_{1}$ and the base eigenvalue $\gl_{0}$ of the hyperbolic Laplacian;
 see \secref{sec:unif}.} And
even more surprisingly,
(2) orbits in 
such
towers
are {\it not} uniformly distributed among different cosets! The uniformity in cosets, were it true, would have 
allowed the application of an
Affine Sieve in this
archimedean
ordering (see, e.g. \cite{Kontorovich2014}); our observation shows that the 
Affine Sieve procedure cannot be  applied directly here, as the standard sieve axioms are not satisfied.\footnote{Of course one can instead order by wordlength in $\G$, as is done in \cite{BourgainGamburdSarnak2010},
to restore equidistribution 
and apply the Affine Sieve.}

\subsection{Statements of the Theorems}\label{sec:intro1}
\

\begin{figure}
        \begin{subfigure}[t]{0.4\textwidth}
                \centering
		\includegraphics[width=\textwidth]{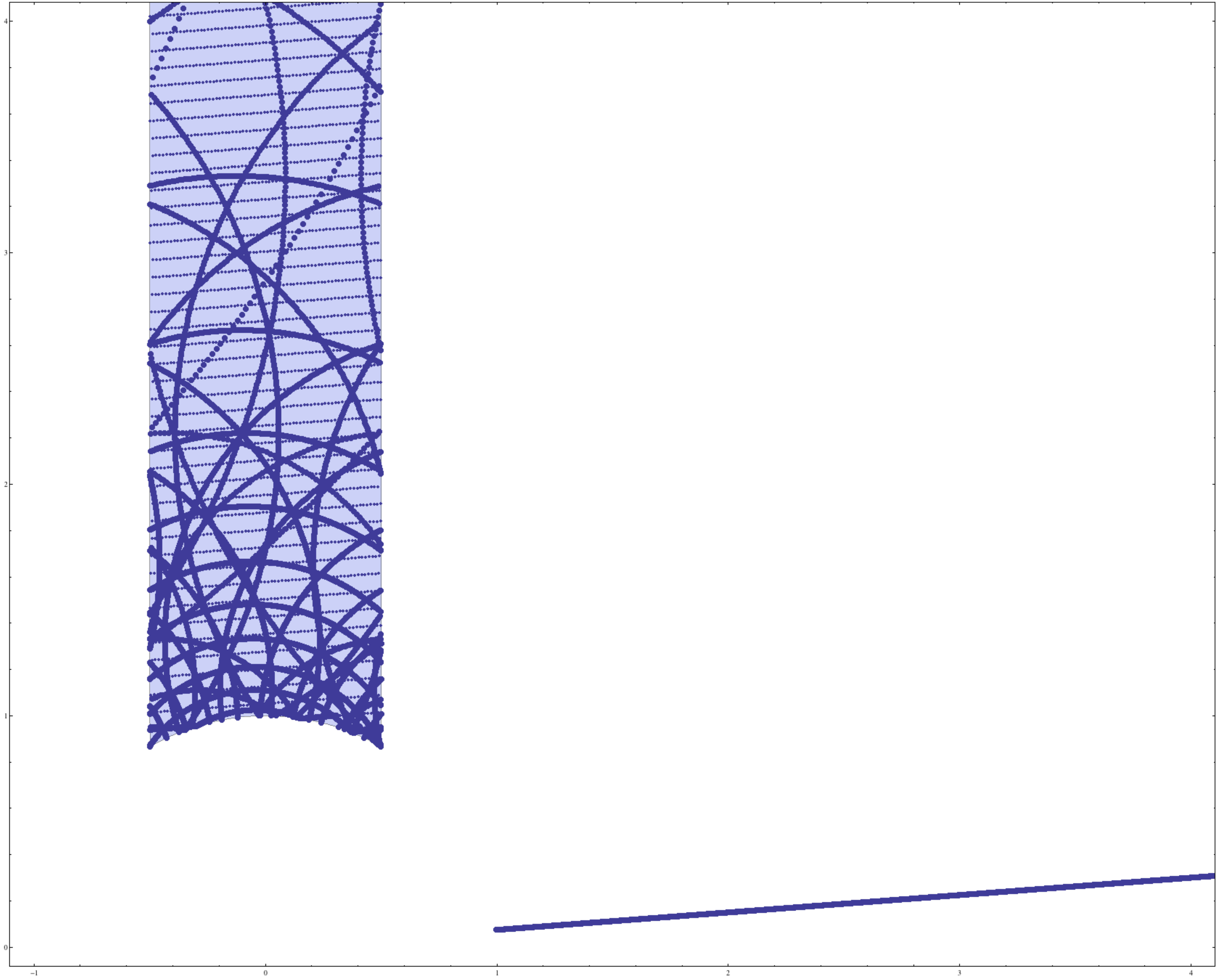}
                \caption{Lattice case: $\G=\PSL_{2}(\Z)$}
                \label{fig:shear}
        \end{subfigure}%
\qquad
        \begin{subfigure}[t]{0.5\textwidth}
                \centering
		\includegraphics[width=\textwidth]{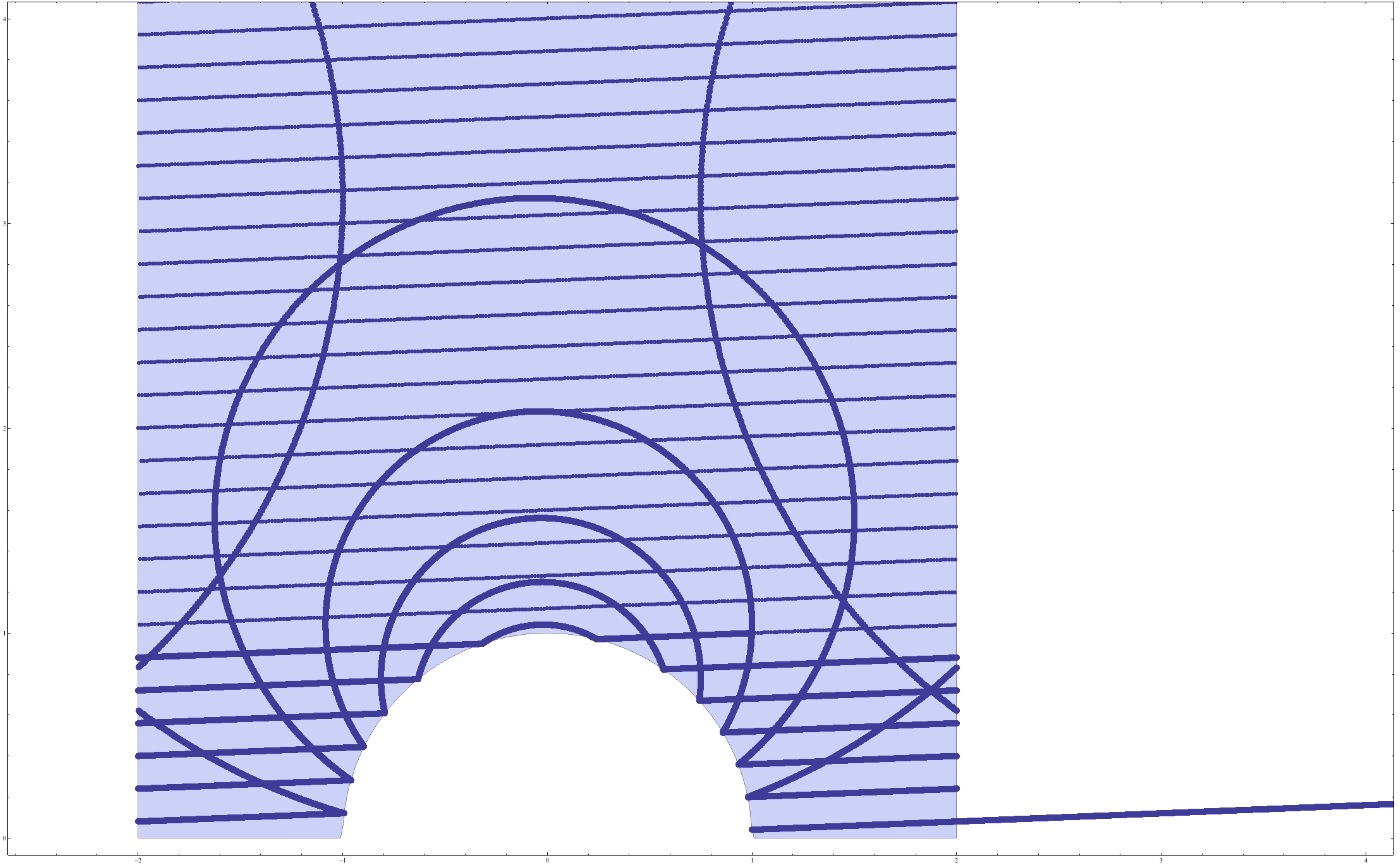}
                \caption{Thin case: $\G=\<\mattwos1401,\mattwos0{-1}10\>$}
                \label{fig:shearThin}
        \end{subfigure}
\caption{A shear of the cuspidal geodesic ray} 
\label{fig:1}
\end{figure}

Our
main
equidistribution
 result 
is the following.
Let $\G$ be a discrete, Zariski-dense,\footnote{Equivalently, non-elementary, that is, not virtually abelian.}
geometrically finite\footnote{For surface groups, 
being geometrically finite is equivalent to being finitely generated.}
subgroup of $G:=\PSL_{2}(\R)$, and assume that the  hyperbolic 
surface
$\G\bk \bH$, which may have finite or infinite volume, has 
at least one 
cusp. In particular, this forces the critical exponent\footnote{%
Roughly speaking, the critical exponent measures the asymptotic growth rate of $\G$; see \secref{sec:hyp}.} $\gd$ of $\G$ to exceed $1/2$; this will be our running assumption throughout.

The base point
$$
\bx_{0}\in T^{1}(\G\bk \bH)\cong\G\bk G
$$
 in the unit tangent bundle 
 determines
 the visual (under the forward geodesic flow) limit point $\fa$ on the boundary $\G\bk\dd\bH$.
 We call the point $\bx_{0}$, 
 as well as
 its 
 forward
 geodesic ray, $\bx_{0}\cdot A^{+}$, {\bf cuspidal} if $\fa$ is a cusp of $\G$; here
 $$
 A^{+}\ := \ \{\mattwos{a}{}{}{a^{-1}}:a>1\}.
 $$
(Note that we make no demands on the  negative geodesic flow from $\bx_{0}$.)
Given such a ray, 
we define its {\bf shear} (the ray is no longer geodesic), at time $T\in\R$,  by:
$$
\bx_{0}\cdot A^{+}\cdot\fs_{T}\ \subset \ \G\bk G
,
$$
where
\be\label{eq:sTdef}
\fs_{T}\ := \ a_{1\over\sqrt{T^{2}+1}}
n_{T},
\qquad\qquad
a_{y} \ = \ \mattwos {\sqrt y}{}{}{1/\sqrt y},
\quad
n_{x} \ =  \ \mattwos 1x{}1.
\ee

For example, if $\G=\SL_{2}(\Z)$, then the base point $\bx_{0}=(i,\uparrow)$ has visual limit point $\fa=\infty$, and hence is cuspidal. The shear at time $T$ of the forward ray from $\bx_{0}$, projected to $\G\bk\bH$, is then simply the Euclidean ray $\{re^{i\gt}\}_{r>1}$, where $\cot\gt=
T$.
See \figref{fig:shear} for an illustration of this ray and its projection mod $\PSL_{2}(\Z)$. Similarly, \figref{fig:shearThin} gives the same picture but for a thin group $\G$.

We are interested in the behavior of such shears as $T\to\infty$ (and similarly for $T\to-\infty$). To this end, define the measure $\mu_{T}$ on a smooth, compactly supported observable $\Psi\in C^{\infty}_{0}(\G\bk G)$ by
\be\label{eq:muTdef}
\mu_{T}(\Psi) \ : = \
\int_{a\in A^{+}}\Psi(\bx_{0}\cdot a \cdot \fs_{T}) da
 \ = \
\int_{1}^{\infty}\Psi(\bx_{0}\cdot a_{y} \cdot \fs_{T}) {dy\over y}
  .
\ee

A slight simplification of our main result (see \thmref{thm:equi}) is the following

\begin{thm}\label{thm:main1}\

{\bf Lattice Case:} Assume that the quotient $\G\bk G$ has finite volume. 
Then
  there exists an $\eta>0$, depending on
   the spectral gap for $\G$,
   so that
\be\label{eq:muFin1}
\mu_{T}(\Psi)
\ = \
\log T 
\cdot
\mu_{\G\bk G}(\Psi)
+
\mu_{\widetilde{Eis}}(\Psi)
+
O_{\Psi}
\left(
T
^{-\eta}\right)
,
\ee
 as $
T\to\infty$.
Here $\mu_{\G\bk G}$ is
the probability
 Haar measure and  $\mu_{\widetilde{Eis}}$ is a certain ``regularized
 Eisenstein''
distribution (see \rmkref{rmk:Etil} below).

{\bf Thin Case:} Assume  that $\G$ is thin, that is, the quotient $\G\bk G$ has infinite volume. Then there exists an $\eta>0$, depending only on
 the critical exponent $\gd$ of $\G$, and \emph{not} on its spectral gap (!), so that
\be\label{eq:muInf1}
\mu_{T}(\Psi)
\ = \
\mu_{{Eis}}(\Psi)
+
O_{\Psi}
\left(
T
^{-\eta}\right)
,
\ee
 as $
 T\to\infty$.
Here $\mu_{Eis}$ is
an 
(un-regularized) Eisenstein distribution.
\end{thm}

Some comments are in order.

\begin{rmk}\label{rmk:tangent}
For simplicity, we have stated \thmref{thm:main1} for compactly supported test functions $\Psi$, but our method applies just as well to a larger class of square-integrable functions with
at least polynomial decay in the cusp $\fa$ (to ensure convergence of $\mu_{T}(\Psi)$); see \secref{sec:mainPf}.
\end{rmk}

\begin{rmk}\label{rmk:optimize}
Throughout we make no attempt to optimize the various error exponents $\eta$, as can surely be done with a modicum of effort; our point is to illustrate a soft method which is powerful enough to obtain new results with power savings errors.
\end{rmk}

\begin{rmk}\label{rmk:Etil}
Let us
 make the 
Eisenstein distributions
arising in
\eqref{eq:muFin1}--\eqref{eq:muInf1}
 less mysterious.
These distributions are actually measures when $\Psi$ is  right $K$-invariant; we restrict attention to this case to simplify the discussion below.
 First assume that $\G$ is a lattice, that $\bx_{0}=(i,\uparrow)$ with $\fa=\infty$ a cusp of $\G$ of width $1$, and 
let $\G_{\infty}
=\mattwos1\Z{}1
$ be the isotropy group of $\fa$ in $\G$. Then one has the standard Eisenstein series
$$
E(z,s)
\ := \
\sum_{\g\in\G_{\infty}\bk\G}
\Im(\g z)^{s}
,
\qquad\qquad(\Re(s)>1),
$$
which is well known
\cite{Selberg1956}
to have meromorphic continuation and a simple pole at $s=1$ with residue $\vol(\G\bk\bH)^{-1}$.
Thus
 the function
\be\label{eq:tilEs}
\widetilde E(z,s) \ : = \
E(z,s)
 \ - \
 {1\over \vol(\G\bk\bH) \ (s-1)}
\ee
is regular at $s=1$;
for example, when $\G=\PSL_{2}(\Z)$, 
we have (see, e.g., \cite[(22.42), (22.63)--(22.69)]{IwaniecKowalski})
\be\label{eq:EisIs}
\widetilde E(z,1)
\ = \
{3\over \pi}
\left(
2\g
-2{\gz'\over\gz}(2)
-\log (4y|\eta(z)|^{4})
\right)
,\qquad
\ee
where $\g=0.577\cdots$ is Euler's constant, $\gz(s)$ is the Riemann zeta function, and $\eta(z)$ is the Dedekind eta function.
Then the measure $\mu_{\widetilde{Eis}}$ is simply given by:
\be\label{eq:tilEisIs}
\mu_{\widetilde{Eis}}(\Psi)
\ = \
\<\Psi,
\widetilde E(\cdot,1)
\>_{\G\bk\bH}
.
\ee
Note that $\log|\eta(z)|$ grows like $y$ in the cusp, so $\mu_{\widetilde{Eis}}$ is also a non-finite measure.
See \eqref{eq:muEis} for the definition when $\Psi$ is not $K$-finite.

In the thin case of \eqref{eq:muInf1}, the Eisenstein series is itself regular at $s=1$, that is, we can simply take $\widetilde{E}(z,s)=E(z,s)$; the spectral contribution is then {\it all} of lower order, so the power savings obtained in \eqref{eq:muInf1} is {\it independent} of any knowledge of a spectral gap for $\G$. 
\end{rmk}

\begin{rmk}\label{rmk:DRS}
The first factor $\log
T$ on the right side of
\eqref{eq:muFin1}
is a manifestation of the logarithmic divergence of the measure $\mu_{T}$.
%
In the lattice case, 
a statement of the form
\be\label{eq:o1}
\mu_{T}(\Psi)
\
=
\
polynomial
(\log T) 
\cdot\mu_{\G\bk G
}(\Psi)
\cdot
\bigg(
1
\ +\
o(
1
)
\bigg)
\ee
was
suggested
 in work of Duke-Rudnick-Sarnak \cite[see below (1.4)]{DukeRudnickSarnak1993}.
Recently,
Oh-Shah
\cite{OhShah2014}
used a purely 
dynamical
method to
prove
(a variant of) 
\eqref{eq:o1}
 with
a  log-
savings
 rate, that is, with $o(
1 )$ replaced by  $O(
1/\log T
 )
 $.
With such a rate it is of course impossible to see the second-order main term (that is, the regularized Eisenstein distribution),
and this identification will be key to some of our applications below.
Moreover, 
it is hard to imagine
how a quantity like \eqref{eq:EisIs} can 
be captured
using
only
dynamics;
our
approach 
 is quite different.
\end{rmk}

Before 
discussing
the proof of \thmref{thm:main1}, we first
give
some applications.

\subsection{Application 1: Weighted Second Moments of $\GL(2)$ $L$-functions}\

Integrals like $\mu_{T}(\Psi)$
arise naturally in
Sarnak's
approach
 (``changing the test vector'') \cite{Sarnak1985a} for second moments of $L$-functions
 (see also, e.g., \cite{
Good1986, Venkatesh2010, MichelVenkatesh2010, DiaconuGarrett2010, DiaconuGarrettGoldfeld2012}).
%
We illustrate the 
method
in the simplest case of
a weight-$k$ holomorphic Hecke cusp form $f$
on $\PSL_{2}(\Z)$, 
though the method works just as well for any $\GL(2)$ automorphic representation.

Write the Fourier expansion of $f$  as
$$
f(z)=\sum_{n\ge1}a_{f}(n)e(nz),
$$
where $a_{f}(1)=1$ and the coefficients $a_{f}(n)$ are multiplicative, satisfying Hecke relations, and the Ramanujan bound  $|a_{f}(p)|\le 2p^{(k-1)/2}$   \cite{Deligne1974}.
The standard $L$-function of $f$,
$$
L(f,s)
\ := \
\sum_{n\ge1}{a_{f}(n)\over n^{s+(k-1)/2}}
,
$$
converges for $\Re(s)>1$, has analytic continuation to $\C$, and a functional equation sending $s\mapsto 1-s$. The Rankin-Selberg 
$L$-function factors (see \cite[(13.1)]{Iwaniec1997book}) as
$$
L(f\otimes\bar f,s)
\ := \
\sum_{n\ge1}{|a_{f}(n)|^{2}\over n^{s+(k-1)}}
 \ = \ {\gz(s)\over\gz(2s)}L(\sym^{2}f,s)
,
$$
where $L(\sym^{2}f,s)$ is the symmetric square $L$-function, known \cite{GelbartJacquet1978} to be the automorphic $L$-function of a cohomological form on $\GL(3)$.

A shear of the 
standard
Hecke integral (already arising implicitly in classical work of Titchmarsch \cite[Chap. VII]{Titchmarsh1951}
) is the following calculation:
\be\label{eq:Hecke}
\int_{0}^{\infty}
f\left(Ty+{iy}\right)
y^{s+(k-1)/2}
{dy\over y}
=
L(f,s)\cW_{k}(s,T)
,
\ee
where
\be\label{eq:cWk}
\cW_{k}(s,T):=
(2\pi)^{-(s+(k-1)/2)}\
\G(s+\tfrac{k-1}2)\
(1-iT)^{-(s+(k-1)/2)}
.
\ee
Applying Parseval to \eqref{eq:Hecke} gives (for $s=1/2+it$)
\be\label{eq:Hecke2}
\int_{0}^{\infty}
\left|
f\left(Ty+{iy}\right)
\right|^{2}
y^{k}
{dy\over y}
=
\frac1{2\pi}
\int_{\R}
|L(f,\tfrac12+it)|^{2}|\cW_{k}(\tfrac12+it,T)|^{2}
dt
.
\ee

\begin{figure}
		\includegraphics[width=.6\textwidth]{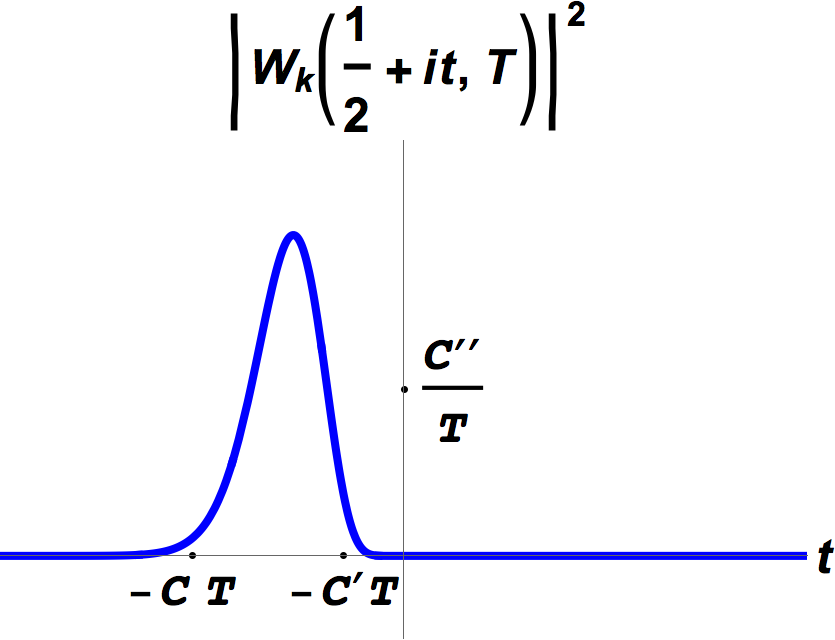}
\caption{A sample graph of the smoothed archimedean weight $|\cW_{k}(\tfrac12+it,T)|^{2}$} 
\label{fig:Wk}
\end{figure}

A calculation with Stirling's formula
(or see  \figref{fig:Wk}) shows that
$|\cW_{k}(\tfrac12+it,T)|^{2}
$
has rapid decay as soon as $|t|>T^{1+\vep}$, and is of size roughly $1/T$ in the bulk.
Thus the quantity on the right side of \eqref{eq:Hecke2}
behaves like a smoothed second moment of $L(f,s)$ on the critical line.
Applying 
\thmref{thm:main1}
with $\Psi=|f|^{2}y^{k}$ (and a little more work, see \secref{sec:thmMom}) gives
the following.
\begin{thm}\label{thm:secMom}
With notation as above, 
there is an $\eta>0$ so that
\bea
\label{eq:thmMom}
&&
\hskip-.5in
\frac1{2\pi}
\int_{\R}
|L(f,\tfrac12+it)|^{2}|\cW_{k}(\tfrac12+it,T)|^{2}
dt
\\
\nonumber
&&
=\
2{\|f 
\|^{2}\over \vol(\G\bk\bH)}
\left(
\log T
\ + \
{\gL'\over \gL}(\sym^{2}f,1)
+
\g
-2
{\gz'\over \gz}(2)
\right)
\ + \
O_{f}(T^{-\eta})
,
\eea
as $T\to\infty$. Here $\g$ is again Euler's constant, $\|f\|$ is the Petersson norm, and $\gL'/\gL$ is the logarithmic derivative of the completed symmetric-square $L$-function,
\be\label{eq:sym2f}
\gL(\sym^{2}f,s)\ = \ (4\pi)^{-(s+k-1)}\,\G(s+k-1)\, L(\sym^{2}f,s). 
\ee
\end{thm}

\begin{rmk}
One can chase the various exponents in our proof to see that \eqref{eq:thmMom} holds with $\eta=1/14-\vep$. Again, we are striving for simplicity of the method and not optimal exponents, see \rmkref{rmk:optimize}.
In fact, a straightforward refinement of the proof of \propref{prop:decayFC} (using explicit spectral expansions in place of soft ergodic arguments) gives $\eta=16/39-\vep$ 
on quoting the best-known bounds \cite{KimSarnak2003} towards the Ramanujan conjectures,
 and $\eta=1/2-\vep$ conditionally.
 So in a sense,
 the proof of \thmref{thm:secMom} is ``sharp,''
 as
  there is no ``loss'' in the rate from a best-possible one.
\end{rmk}

\begin{rmk}
On comparing
the lower order terms on 
the right hand side of \eqref{eq:thmMom} 
with the secondary term 
 in  \eqref{eq:muFin1},
and using \eqref{eq:EisIs}
and \eqref{eq:tilEisIs},
one
derives 
a
 Kronecker-type limit formula, in the form:
\be\label{eq:logEta}
{
\<
\log (4y|\eta(z)|^{4})
,|f|^{2}y^{k}\>
\over
 \|f\|^{2}}
\ =\
\g
\ - \
{\gL' \over \gL}(\sym^{2}f,1)
.
\ee
This 
identity
is surely classical, 
though
we were not able to locate a precise reference. 
\end{rmk}

\subsection{Application 2: Archimedean Counting for Orbits on Affine Quadrics}\

Another standard
context where integrals like $\mu_T 
(\Psi)$ arise
naturally
 is in
 the execution of certain
  Margulis/Duke-Rudnick-Sarnak/Eskin-McMullen type arguments
\cite{Margulis2004, DukeRudnickSarnak1993,
EskinMcMullen1993}
for
 counting 
 discrete orbits
 on quadrics
in  archimedean balls.
%
%
The setting is as follows.

Let $Q$ be a real ternary indefinite 
quadratic form
(e.g.,
 $Q(\bx)=x^{2}+y^{2}-z^{2}$),
fix
 $d\in\R$, 
and denote by $V=V_{Q,d}$ the affine quadric
\be\label{eq:VQd}
V\ : \ Q=d.
\ee
The real points $V(\R)$ enjoy
 an action by $G=\SO_{Q}^{\circ}(\R)$, the connected component of the identity in the real special orthogonal group preserving $Q$. Let
$\G<G$ be a discrete, Zariski dense, geometrically finite subgroup of $G$, and
assume, as throughout, that
the critical exponent $\gd$ of $\G$ exceeds $1/2$.
Fix a base point $\bx_{0}\in V(\R)$, 
subject to
 the orbit
$$
\cO\ :=\ \bx_{0}\cdot\G\ \subset\ \R^{3}
$$
being discrete.

For a fixed archimedean norm $\|\cdot\|$ on $\R^{3}$, let
$$
B_{T}=\{\bx\in\R^{3}:\|\bx\|<T\}
$$
be the norm ball of radius $T$.
A very classical and well-studied problem is to
give an effective (that is, with power savings error)
estimate for
  $$
\cN_{\cO
}(T)
\ :=\
\left|
\cO 
\ \cap \ B_{T}
\right|
.
  $$

Despite the vast attention this problem has received over the years, there remained exactly two lacunary cases in which hitherto resisted
solution; see
 \secref{sec:taxonomy} 
 and \tabref{tab:1}
 below
for a detailed taxonomy of the situation.
Equipped with \thmref{thm:main1}, we can now resolve the outstanding cases. 

\begin{thm}\label{thm:count}
There exist constants $C_{1}, C_{2},$ and $\eta>0$ so that the following holds:
\begin{itemize}
\item If $\G$ is a lattice in $G$, then
\be\label{eq:countLat}
\cN_{\cO}(T) \ = \
\big(C_{1}T\log T + C_{2} T\big)\big(1 +O(T^{-\eta})\big)
.
\ee
\item If $\G$ is
thin, then
\be\label{eq:countThin}
\cN_{\cO}(T) \ = \
\big(C_{1}T + C_{2} T^{\gd}\big)\big(1 +O(T^{-\eta})\big)
.
\ee
\end{itemize}
\end{thm}

Some comments are in order:
\begin{rmk}\label{rmk:DRSEM}\

\begin{enumerate}
\item
All of the previously
resolved
cases of this problem (in the above generality) were such that the first term did not appear, that is, $C_{1}=0$ (whence $C_{2}>0$). Our new contributions are to the cases with $C_{1}>0$, which arise exactly when the stabilizer of $\bx_{0}$ in $G$ is
diagonalizable but
 ``cuspidal'' in $\G\bk G$; see 
 \secref{sec:taxonomy}
 below. 

\item
If it happens that
 $\cO$ is not just some arbitrary real discrete orbit but is actually the full integer quadric $V(\Z)$ (assuming of course that the quadratic form $Q$ is rational and that $V(\Z)$ is non-empty),
 then
 one has many more tools available  to approach the counting problem for $\cN_{\cO}(T)$. Specifically, one can use, e.g.,  classical methods of exponential sums (see \cite{Hooley1963}), or half-integral weight automorphic forms, Poincar\'e series and  shifted convolutions \cite{Sarnak1984,  Blomer2008, TemplierTsimerman2013}, or  
multiple Dirichlet series  \cite{HulseKiralKuanLim2013}. These give, when $\G$ is a lattice and $C_{1}>0$, an estimate for $\cN_{\cO}(T)$ of the same strength as \eqref{eq:countLat}, that is, with a secondary main term and  a power savings error. 
Such
tools do not seem to apply to the general orbit counting problem.

\item
In the thin case with $C_{1}>0$, the second term $C_{2}T^{\gd}$ is  swamped by the error, and should not be confused with a lower order ``main'' term.
We would like to acknowledge here that Nimish Shah suggested to us that the main term in this setting is of order $ T$ rather than $T^{\gd}$; see also \cite[Remark 1.7]{OhShah2013}.

\item
For the ``new'' cases with $C_{1}>0$, the exponent $\eta$ depends on the same quantities as  in \thmref{thm:main1}; that is,
$\eta$ depends on the spectral gap in the lattice case, and
 only on the critical exponent  in the thin case.

\item
The constants $C_{1}, C_{2}$ can be readily determined explicitly
in terms of volumes, special values of (possibly regularized) Eisenstein series, and Patterson-Sullivan measures.

\item
A consequence of Oh-Shah's result
discussed in \rmkref{rmk:DRS}
 gives,
  in the lattice case, 
the   estimate
  \eqref{eq:countLat},
but with
$O(T^{-\eta})$ replaced by
 the
weaker error rate  $O(1/(\log T)^{\eta})$ for some small $\eta>0$. 
This of course only identifies the first main term $C_{1}T\log T$, as the secondary term $C_{2}T$ is swallowed by the error.

\end{enumerate}
\end{rmk}


\subsubsection{Taxonomy}\label{sec:taxonomy}\

To explain the lacunary cases settled by \thmref{thm:count}, we  begin by
passing from $\SO_{Q}^{\circ}(\R)$ to its spin cover $\PSL_{2}(\R)\cong T^{1}(\bH).$
Abusing notation, we continue to write $G$ and $\G$
for their pre-images in $\PSL_{2}(\R)$.

Let $H$ be the stabilizer of $\bx_{0}$ in $G$,
$$
H \ := \ \{h\in G : \bx_{0}\cdot h=\bx_{0}\},
$$
 and let 
 $$
 \G_{H}:=\G\cap H.
 $$ 
With $\bx_{0}$ fixed, the norm $\|\cdot\|$ on $\R^{3}$ induces a left-$H$ invariant norm $\vertiii{\cdot}$ on $
G$ given by 
\be\label{eq:normiii}
\vertiii{g}=\|\bx_{0}g\|.
\ee 
We further abuse notation, writing $B_{T}$ for the left-$H$-invariant norm-$T$ ball in $ G$, that is, $B_{T}\subset H\bk G$.
Then it is easy to see that
$$
\cN_{\cO 
}(T)
\ = \
|
\G_{H}\bk \G 
\ \cap \ B_{T}
|
.
$$
%

\begin{figure}
        \begin{subfigure}[t]{0.3\textwidth}
                \centering
		\includegraphics[height=1.6in]{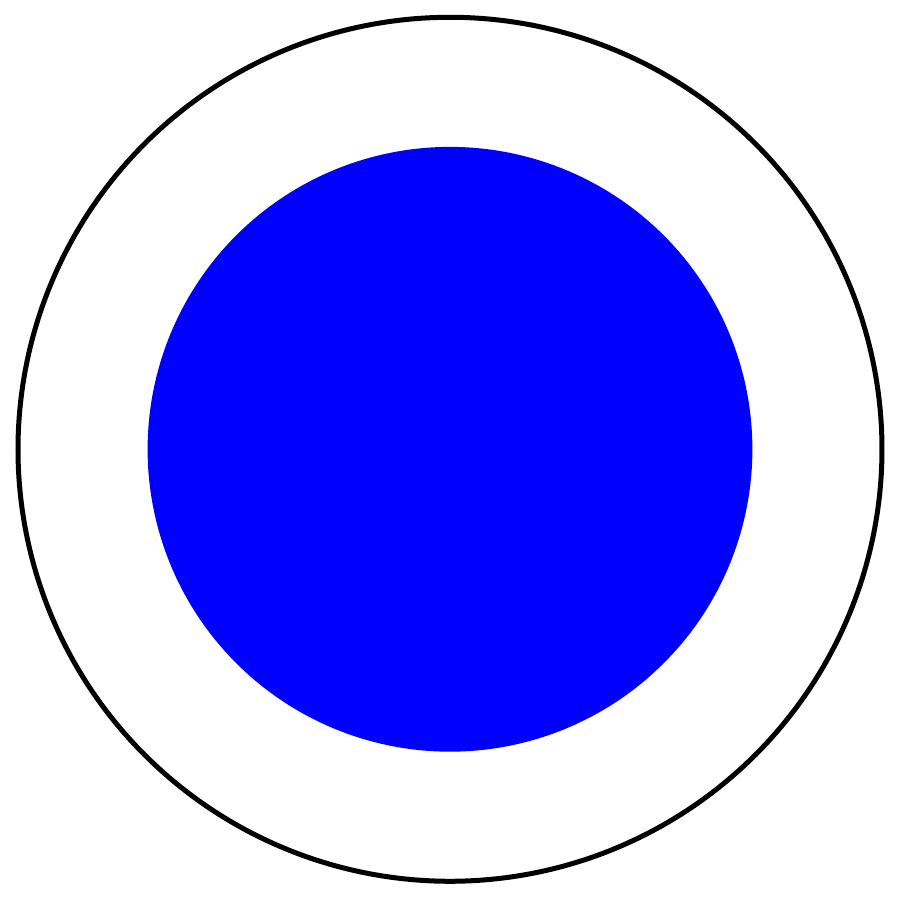}
                \caption{Case $H\cong K$.}
\label{fig:BTK}        
\end{subfigure}%
\quad
        \begin{subfigure}[t]{0.3\textwidth}
                \centering
		\includegraphics[height=1.6in]{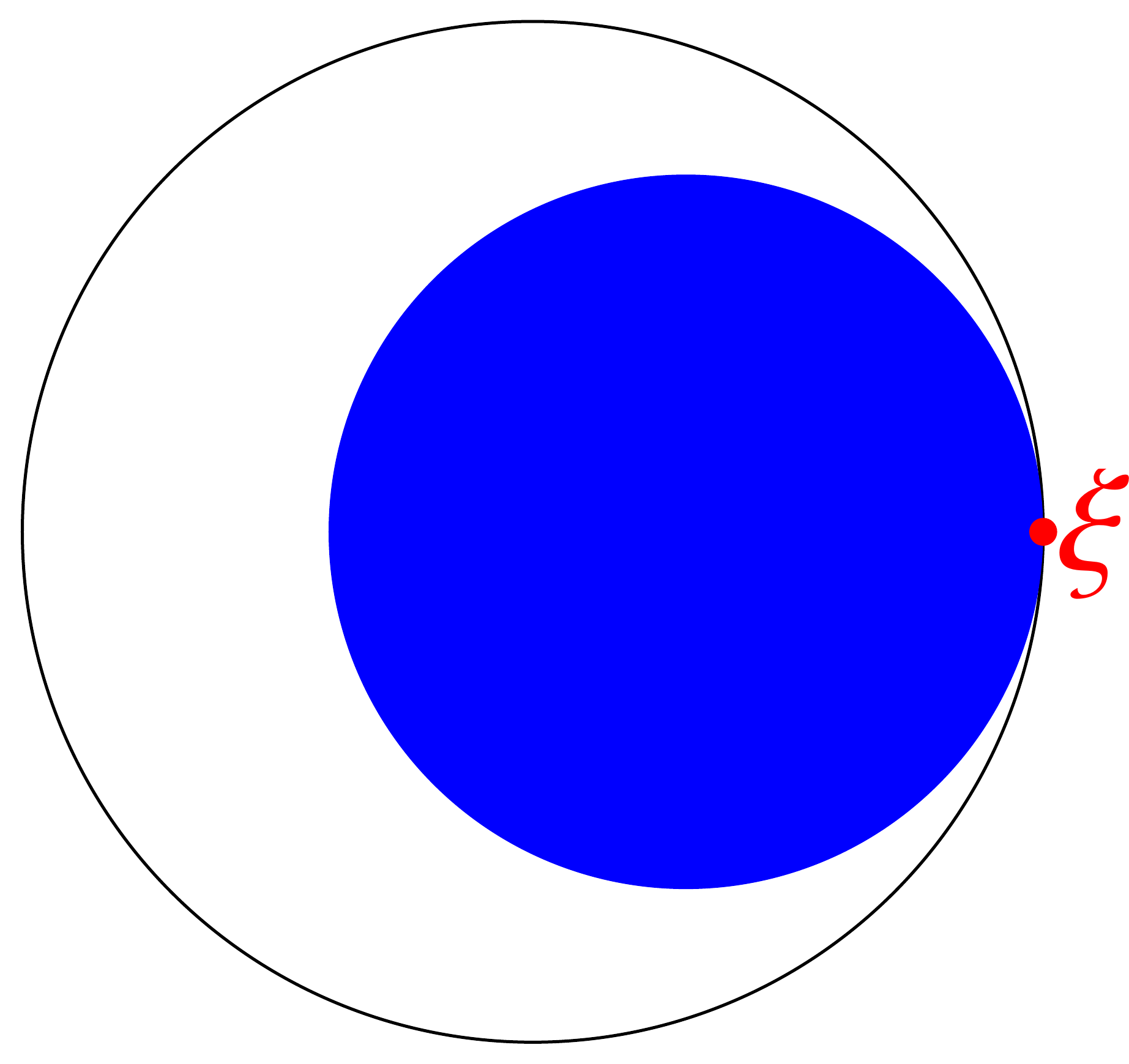}
                \caption{Case $H\cong N$.}
\label{fig:BTN}        
        \end{subfigure}%
\quad
        \begin{subfigure}[t]{0.3\textwidth}
                \centering
		\includegraphics[height=1.6in]{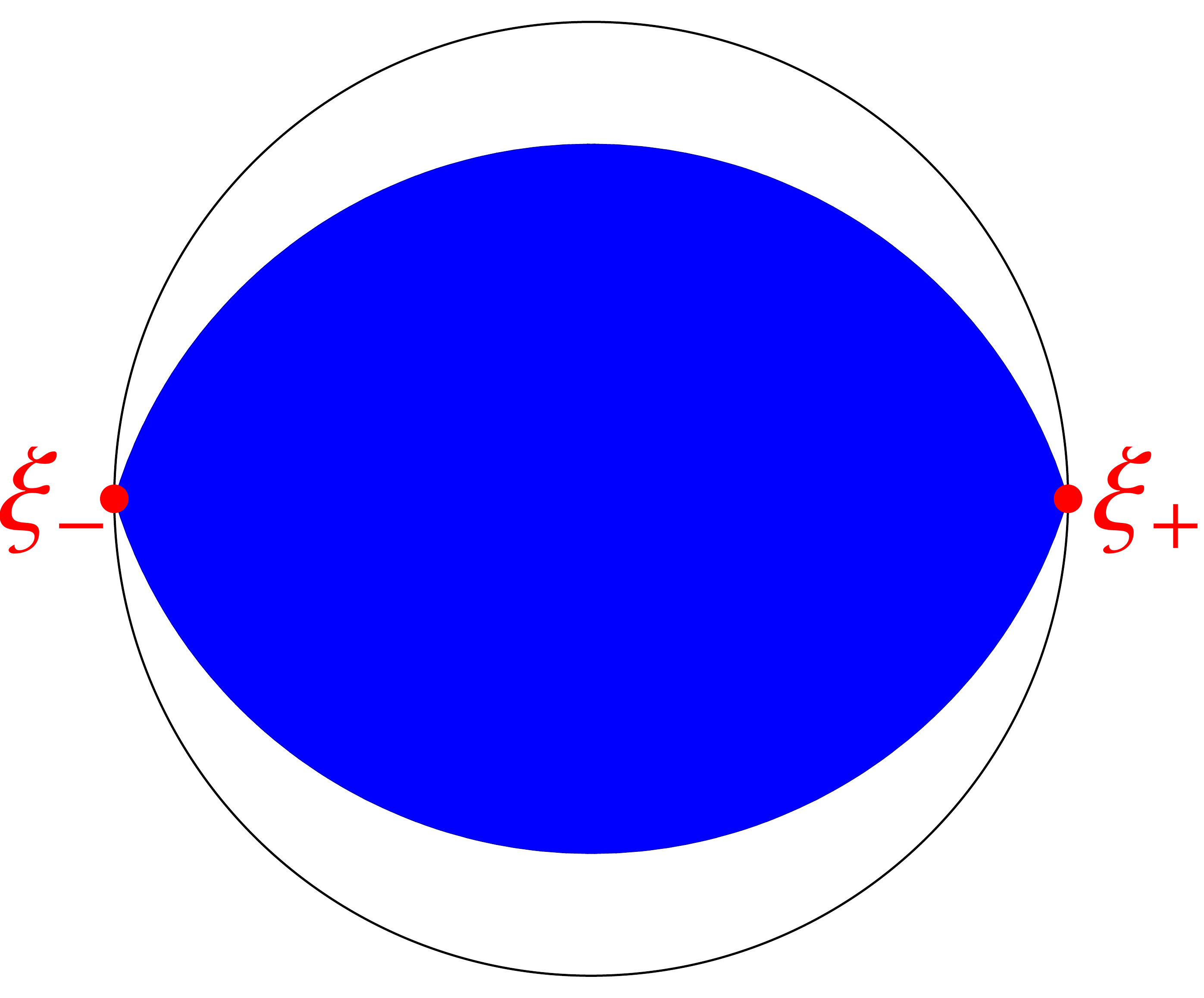}
                \caption{Case $H\cong A$.}
\label{fig:BTA}        
        \end{subfigure}%
\caption{The region $B_{T}$ as a subset of the hyperbolic disk $\bD$.} 
\label{fig:BT}
\end{figure}

To investigate
 this counting problem
 more 
 precisely,
 we 
 illustrate the geometry of $B_{T}$, 
which is
 determined by whether the stabilizer $H$ is conjugate to
 the groups 
 $K$, $N$, or $A$. 
That is, $H$ is either a maximal compact, 
a unipotent subgroup, or
diagonalizable (over $\R$),
and this corresponds to whether
the real quadric $V(\R)$ 
is 
a two-sheeted hyperboloid, a  cone, 
or a
one-sheeted hyperboloid, 
respectively.
To visualize $B_{T}$ as a left-$H$-invariant subset of $G$,
we project to the base space $\bH$ (or alternatively,
 assume that the norm $\vertiii{\cdot}$ is right-$K$-invariant), so that $B_{T}
 $ can be viewed as an $H$-invariant subset of the hyperbolic disk $\bD\cong G/K$. Then $B_{T}$ is 
illustrated
in \figref{fig:BT} in the three cases.
Note that $B_{T}$ has zero, one, or two limit points (denoted $\xi$ or $\xi_{\pm}$) on the 
 boundary $\dd\bD$, corresponding to whether  $H\cong K$, $H\cong N$,
or $H\cong A$,
respectively.

The asymptotic
counting
 analysis 
depends in a fundamental way not only on whether $\G$ is a lattice in $G$, but also on whether
\be\label{eq:ass}
\text{$\G_{H}$ is a lattice in $H$}.
\ee

\begin{lem}\label{lem:cO}
If \eqref{eq:ass} does not hold, then  the discreteness of $\cO$ is equivalent to the endpoints $\xi$ or $\xi_{\pm}$ of $H$ not being radial limit points\footnote{Recall that the limit set, $\gL$, of $\G$ decomposes disjointly into cusps (i.e., parabolic fixed points) and radial limit points (also called ``points of approximation''); the complement $\dd\bH\setminus\gL$ is called the free boundary (which is empty if $\G$ is a lattice). See \secref{sec:hyp}.} for $\G$. 
\end{lem}

This follows from a simple topological argument; we omit the proof. We decompose the analysis according to whether $\G$ is a lattice or thin in $G$.

{\bf Case I:
 $\G$ is a lattice.} \

Assuming that $\G$ is a lattice in $G$, and also 
demanding
 that
\eqref{eq:ass} holds,
Duke-Rudnick-Sarnak \cite{DukeRudnickSarnak1993} 
and Eskin-McMullen \cite{EskinMcMullen1993}
(see also \cite{Margulis2004}) showed (in much greater generality than considered here) that there is some $\eta>0$
with
\be\label{eq:DRS}
\cN_{\cO}(T)
\ =\
{\vol_{H}(\G_{H}\bk H)\over \vol_{G}(\G\bk G)}
\vol_{H\bk G}(B_{T})
\Bigg(
1+O(T^{-\eta})
\Bigg)
,
\ee
as $T\to\infty$.
Here the volumes are taken to be compatible with choices of Haar measure on $G$, 
 $H$, and $H\bk G$.
Note that $\vol_{H\bk G}(B_{T})$ is of order $T$, so there is no logarithmic divergence in \eqref{eq:DRS}, that is, $C_{1}=0$
and $C_{2}>0$
 in \eqref{eq:countLat};
see also \rmkref{rmk:DRSEM}(1).

With \figref{fig:BT} and \lemref{lem:cO} in mind, we 
analyze separately the possible conjugacy classes of $H$.
\begin{itemize}
\item First if $H\cong K$ is compact, then 
\eqref{eq:ass} clearly holds automatically.
In this case, the counting result 
\eqref{eq:DRS} 
 corresponds to counting in norm balls of $G$, which 
dates back
to Delsarte \cite{Delsarte1942}, Huber \cite{Huber1956}, and Selberg \cite{Selberg1956}.
\item If $H\cong N$ is unipotent, then by \lemref{lem:cO}, the boundary point $\xi$ of $H$ must be a cusp. 
That is,
 $\G_{H}\bk H$ is a closed horocycle, so
 $\G_{H}$ is a lattice in $H$, and
 \eqref{eq:ass} is again automatically satisfied. 
In this case,
the count takes place in a strip $\G_{H}\bk G$, and
 the equidistribution of low-lying closed horocycles \cite{Margulis2004, Zagier1981, Sarnak1981} can be used to establish \eqref{eq:DRS}.
\item
Lastly, if $H\cong A$ is diagonalizable (over $\R$), then 
\lemref{lem:cO}
forces one of two settings. 
Either 

\subitem
$(i)$: $\G_{H}$ is a lattice in $H$, whence $\G_{H}\bk H$ corresponds to a closed geodesic on $\G\bk G$. Then \eqref{eq:ass} holds, so \eqref{eq:DRS} follows from \cite{DukeRudnickSarnak1993}. Or

\subitem
$(ii)$:
$\G_{H}$ is finite, but
both limit points $\xi_{+}$ and $\xi_{-}$ of $H$ (see  \figref{fig:BTA}) are cusps of $\G$.
Here we are in  the diagonalizable but ``cuspidal'' setting referred to in \rmkref{rmk:DRSEM}(1).
 This 
 case
is the only one (when $\G$ is a lattice in $G$) not satisfying \eqref{eq:ass} despite the discreteness of the orbit $\cO$; it is precisely
 the new case 
  settled by \thmref{thm:count}.  
\end{itemize}

\

{\bf Case II: $\G$ is thin.}\

In this setting, we again decompose the problem of estimating $\cN_{\cO}(T)$ into separate cases, depending
 on the conjugacy class of $H$, and on whether condition \eqref{eq:ass} holds (there are now more 
 situations
 in which $\cO$ is discrete but \eqref{eq:ass} can fail). 
 \begin{itemize}
 \item
When $H\cong K$ is compact, the condition \eqref{eq:ass} is again automatically satisfied, and in this case Lax-Phillips \cite{LaxPhillips1982} showed that
\be\label{eq:LaxP}
\cN_{\cO}(T)=C_{2} \, T^{\gd}\Bigg(
1+O(T^{-\eta})
\Bigg)
,
\ee
where $\eta>0$ depends on the spectral gap for $\G$.
This corresponds to the case $C_{1}=0$
 in \eqref{eq:countThin}; again,
see \rmkref{rmk:DRSEM}(1).

 \item
If $H\cong N$ is unipotent, the discreteness of $\cO$ forces one of two cases. Either

\subitem
(i) $\G_{H}$ is a lattice in $H$, so \eqref{eq:ass} holds, and $\G_{H}\bk H$ corresponds to a closed horocycle. In this case, the asymptotic formula is   \eqref{eq:LaxP} was shown in the second-named author's thesis \cite{Kontorovich2009}. Or

\subitem
(ii) $\G_{H}$ is trivial, whence \lemref{lem:cO} forces the limit point $\xi$ of $H$ 
to be
in the free boundary, that  is, it is not in the limit set of $\G$. 
The asymptotic here
is
also 
of the form \eqref{eq:LaxP}; 
 see
 \cite{KontorovichOh2012}.

\item
Finally, when $H\cong A$  is diagonalizable,
 there are three separate cases to consider. The discreteness of $\cO$ now implies either

\subitem
(i)
 $\G_{H}$ is a lattice in $H$, again corresponding to a closed geodesic on $\G\bk G$.
Then 
the same asymptotic  \eqref{eq:LaxP} follows from now-standard methods using the equidistribution result of Bourgain-Kontorovich-Sarnak \cite{BourgainKontorovichSarnak2010}. Or

\subitem
(ii)
$\G_{H}$ is  thin in $H$, in which case each of the two endpoints $\xi_{\pm}$ of $H$ is either a cusp of $\G$ or in the free boundary. If 

\subsubitem
(a) both endpoints $\xi_{\pm}$ are in the free boundary, then the methods of \cite{BourgainKontorovichSarnak2010} can again be used to show the same asymptotic \eqref{eq:LaxP}. (The key
is that
 only a finite portion of the sheared geodesic ray interacts with the limit set -- see \cite{KontorovichOh2012} where a similar phenomenon was studied in the case of a unipotent stabilizer.) 
Otherwise,
 
\subsubitem
(b) at least one of $\xi_{\pm}$ is a cusp of $\G$. This is the other %
new
lacunary
 case of \thmref{thm:count}, and is the only
  thin case for which \eqref{eq:countThin} has $C_{1}>0$. Note that if one boundary point is a cusp and the other is in the free boundary, the former has contribution of order $T$, while the latter's contribution, of order $T^{\gd}$, 
is dominated by the former's error term.

\end{itemize}

This concludes our taxonomy. 
To summarize, the following table serves to illustrate that Theorem \ref{thm:count} above completes 
the effective solution to the general counting orbital problem in our context:

\begin{table}[h]
\hspace*{-50 pt}
\begin{tabular}{|c|c|c|c|}
\hline
\backslashbox{$(\G,\G_{H})$}{$H
$ 
}&$N$&$A$&$K$\\
\hline
(lattice, lattice)&\cite{
Zagier1981, Sarnak1981}&\cite{DukeRudnickSarnak1993, EskinMcMullen1993}&\cite{Delsarte1942, Huber1956, Selberg1956}\\
\hline
(lattice, thin)&
\begin{tabular}{c}
impossible by \\
discreteness of $\cO$
\end{tabular}
&
\cellcolor{pink}{
\begin{tabular}{c}
``lacunary'' case \\
settled in \eqref{eq:countLat}
\end{tabular}
}
&
\begin{tabular}{c}
impossible by \\
compactness of $K$
\end{tabular}
\\
\hline
(thin, lattice)&
\cite{Kontorovich2009}
&\cite{BourgainKontorovichSarnak2010}&
\cite{LaxPhillips1982}
\\
\hline
(thin, thin)&
\cite{KontorovichOh2012}
&
$
\twocase{}
{\text{\cite{BourgainKontorovichSarnak2010},}}{if both $\xi_{\pm}\notin\gL$,  }
{\cellcolor{pink}
{
%
%
\hskip-.1in
\begin{tabular}{c}
``lacunary'' case \\
settled in \eqref{eq:countThin},
\end{tabular}
}
}{
otherwise,
}
$
&
\begin{tabular}{c}
impossible by \\
compactness of $K$
\end{tabular}
\\
\hline
\end{tabular}
\caption{The new cases of \thmref{thm:count},  highlighted, are those with $C_{1}>0$.}
\label{tab:1}
\end{table}


\begin{rmk}
When the critical exponent $\gd\le 1/2$, work of Naud \cite{Naud2005}, extending Dolgopyat's methods \cite{Dolgopyat1998}, allows one to conclude, in the cases not excluded by \lemref{lem:cO}, an effective asymptotic of the form \eqref{eq:LaxP}. So the lacunary 
cases
do not occur 
here, 
as there are no cusps (and hence no cuspidal geodesic rays) when $\gd\le1/2$.
\end{rmk}

\begin{rmk}
As pointed out in \cite[p. 917]{OhShah2014} (at least for $\G$ a lattice), the only lacunary cases in the more general setting of $Q$ having signature $(n,m)$ are precisely those of signature $(2,1)$, that is, those considered here; so we have lost no generality in restricting to $\PSL_{2}(\R)$. This is because the stabilizer $H$ is either unipotent, compact, or  fixes a 
form of signature $(n-1,m)$ or $(n,m-1)$. The only
 non-compact
  such 
  not generated by unipotents has signature $(1,1)$, whence $Q$ has signature $(2,1)$. 
\end{rmk}

\subsection{Surprise: Non-equidistribution in Congruence Cosets!}\label{sec:introCosets}\

The most unexpected consequence of \thmref{thm:count} comes from studying the thin case in cosets of congruence towers, as we now describe.
Assume  that
$Q$ is not just a real quadratic form but an integral one, and that
 $\G$ is
 a subgroup of the {\it integral} special orthogonal group
$$
\G\ < \ \SO_{Q}(\Z).
$$
Given an integer $q\ge1$, we can then form the level-$q$ ``congruence'' subgroups $\G(q)<\G$, defined as
$$
\G(q) \ : =  \ \ker(\G\to\SO_{Q}(\Z/q\Z)).
$$
For many applications, one wishes to study the same counting problem as above, with the orbit $\cO=\bx_{0}\cdot\G$ replaced by the congruence orbit $\bx_{0}\cdot\G(q)$, or better yet, by some ``congruence coset'' orbit,
$$
\cO_{q,\vp} \ := \ \bx_{0}\cdot\vp\cdot\G(q),
$$
for a given $\vp\in\G/\G(q)$. Let
$$
\cN_{q,\vp}(T ) \ := \ |\cO_{q,\vp} \ \cap \ B_{T} |
$$
be the corresponding counting function, which we wish to estimate {\it uniformly} with $q$ and $T$ (and $\vp$) varying in some allowable range.

\thmref{thm:count} applies just as well to estimate $\cN_{q,\vp}(T)$, and in all previously studied examples, the asymptotic analysis showed that
\be\label{eq:edcosets}
\cN_{q,\vp}(T) \ \sim \ \frac1{[\G:\G(q)]}\cN_{\cO}(T),\qquad\qquad\qquad (T\to\infty),
\ee
that is,
the asymptotic is
 independent of $\vp$, so
 the orbits are 
 {\it equidistributed} among congruence cosets. (Moreover, \eqref{eq:edcosets} even holds with $q$
allowed to
grow 
sufficiently slowly with $T$.)
This equidistribution is a key input, for example, in executing an Affine Sieve in an archimedean ordering (see, e.g., \cite{Kontorovich2009,  NevoSarnak2010, LiuSarnak2010, BourgainGamburdSarnak2011,  KontorovichOh2012, MohammadiOh2013, BourgainKontorovich2014b, HongKontorovich2014}).
An
analysis of \thmref{thm:count} shows that, for thin orbits, there are cosets
 which  are {\it not} uniformly distributed in archimedean balls!

\begin{prop}\label{prop:noED}
Assume that $\G$ is thin, and that the orbit $\cO$ is
has diagonalizable and cuspidal stabilizer, that is,
 $C_{1}>0$ in \eqref{eq:countThin}. Then the equidistribution \eqref{eq:edcosets} in congruence cosets is {\bf false}. For example, for each fixed $q$, there is some $\vp\in\G/\G(q)$ so that
\be\label{eq:nonEDqT}
\cN_{q,\vp}(T) \ \gg \ \frac 1q\ T,
\ee
while
$$
\frac1{[\G:\G(q)]}\cN_{\cO}(T) \ \ll \ \frac 1{q^{3}} \ T,
$$
as $T\to\infty$. (The implied constants above may depend on $\G$ and $\bx_{0}$, but obviously not on $q$ or $T$.)
\end{prop}

This means that the standard  Affine Sieve procedure {\it cannot} be executed in this ordering. 
(Note that the case considered in \propref{prop:noED}
is precisely  the one omitted in the sieving statement  \cite[Cor 1.19]{MohammadiOh2013}.)



\subsection{Outline of the Proofs and Paper}\

Our proof of \thmref{thm:main1} is surprisingly 
simple, and proceeds in two stages. The first is to show  that $\mu_{T}(\Psi)$ in some sense equidistributes in the ``strip $\G_{\infty}\bk \bH$'', but with respect to $dx\,dy/y$, as opposed to
 Haar measure
 $dx\,dy/y^{2}$
  (for a hint of this, look again at \figref{fig:1}); this fact uses only the decay of the Fourier coefficients of $\Psi$
(
 itself
 a simple consequence of mixing 
 via
 low-lying horocycles; see \propref{prop:decayFC}). Then stage two is to relate this equidistribution to
  Eisenstein series,
  where we
  mimic Sarnak's approach in \cite{Sarnak1981, Zagier1981} to conclude the proof.

The rest of the paper proceeds as follows.
In \secref{sec:bkg}, we set notation and recall basic facts needed through the paper. Then in \secref{sec:mainPf}, we prove the main equidistribution result, \thmref{thm:main1}, and its generalization (\thmref{thm:equi}). 
Then \thmref{thm:secMom} is proved in \secref{sec:thmMom} using the Rankin-Selberg ``unfolding'' technique. Finally, \thmref{thm:count} and \propref{prop:noED} are proved in \secref{sec:counting}.

\subsection{Notation}\

Constants $0<C<\infty$ and $0<\eta<1$ can change from line to line, and $\vep>0$ represents an arbitrarily small quantity. The transpose of a matrix $g$ is written $^{\top}\!g$.
Unless otherwise specified, implied constants depend at most on $\G$, which is treated as fixed. The symbol $\bo_{\{\cdot\}}$ represents the indicator function of the event $\{\cdot\}$.

\subsection*{Acknowledgements}\

The authors would like to express their gratitude to
Jens Marklof, Curt McMullen, and Peter Sarnak for many enlightening comments and suggestions. The second author would especially like to thank
 Valentin Blomer, Farrell Brumley, and Nicolas Templier for many hours of discussion about this problem; see also the related work in \cite{BlomerBrumleyKontorovichTemplier2014}.

\newpage

\section{Preliminaries}\label{sec:bkg}

In this section, we set all notation and basic facts used throughout.

\subsection{Hyperbolic Geometry}\label{sec:hyp}\

Let $\bH:=\{z\in\C:\Im z>0\}$ denote the hyperbolic upper half-plane. At each point $z\in\bH$, and tangent vector $\gz\in T_{z}\bH\cong\C$, the Riemannian structure is
$
\|\gz\|_{z} :=|\gz|/\Im z.
$
The unit tangent bundle $T^{1}\bH$ is then
$$
T^{1}\bH \ := \ \{(z,\gz)\in\bH\times\C:\|\gz\|_{z}=1\}.
$$
The group $G=\PSL_{2}(\R)$ acts on $T^{1}\bH$ via
$$
\mattwo abcd\ : \ (z,\gz)\mapsto \left({az+b\over cz+d},{\gz\over (cz+d)^{2}}\right),
$$
and moreover we can identify $G\cong T^{1}\bH$ under $g\mapsto g(i,\uparrow)$.
We also use the disk model $\bD:=\{z\in\C:|z|<1\},$ identified with $\bH$ under the map
$$
\bH\ni z \ \mapsto \ (z-i)/(z+i)\in\bD.
$$

Let $\G$ be a finitely generated,  Zariski dense,  discrete subgroup of $G$. As above, we identify $T^{1}(\G\bk\bH)\cong \G\bk G$.
 For a fixed base point $\fo\in\bH$, the {\bf critical exponent}
$$
\gd=\gd(\G) \ \in \ [0,1]
$$
of $\G$ is the abscissa of convergence of the Poincar\'e series
$$
\sum_{\g\in\G}\exp(-sd(\g\fo,\fo)), \hskip1in(\Re(s)>\gd).
$$
Here $d(\cdot,\cdot)$ is hyperbolic distance, and $\gd$ does not depend on the choice of $\fo$.
Let $dg$ be a choice of Haar measure on $G$; we call $\G$ a {\bf lattice} if $\G\bk G$ has finite measure, and {\bf thin} otherwise.
This is measured by the critical exponent $\gd$, as $\gd=1$ or $\gd<1$ exactly when $\G$ is a lattice or thin, respectively \cite{Patterson1976};
the Zariski-density of $\G$ implies that $\gd>0$.
The {\bf limit set}
$$
\gL=\gL(\G)\ \subset \ \dd\bH\cong S^{1}\cong\R\sqcup\{\infty\}
$$
of $\G$ is the set of limit points of $\g\fo$, $\g\in\G$; it also does not depend on the choice of $\fo$.
The Hausdorff dimension of $\gL$ is exactly equal to the critical exponent $\gd$ \cite{Patterson1976, Sullivan1984}.
A boundary point $\xi\in\dd\bH$ is a {\bf cusp} of $\G$ if it is the fixed point of  a parabolic element in $\G$; these all lie in the limit set $\gL$, and we let $\gL_{cusp}$ denote the subset of cusps. A limit point $\xi\in\gL$ is called {\bf radial} (or a ``point of approximation'') if there is a sequence $\{\g_{n}\fo\}$, $\g_{n}\in\G$, which stays a bounded distance away from a geodesic ray ending at $\xi$. Let $\gL_{rad}$ denote the set of radial limit points;
it is a basic fact
 \cite{Beardon1983} that the limit set decomposes disjointly into radial and cuspidal points,
 $$
 \gL  \ = \ \gL_{cusp}\sqcup\gL_{rad}
 .
$$
The complement of the limit set in $\dd\bH$ is called the {\bf free boundary} of $\G$,
$$
\cF=\cF(\G) \ : = \ \dd\bH \setminus\gL,
$$
and $\cF=\O$ if and only if $\G$ is a lattice.
We record here the decomposition
\be\label{eq:ddbH}
\dd\bH \ = \ \cF \sqcup \gL_{cusp}\sqcup\gL_{rad}.
\ee
We assume henceforth that $\G$ has at least one cusp, whence its critical exponent exceeds $1/2$,
\be\label{eq:del12}
\gd \ > \ 1/2.
\ee

\subsection{Spectral Theory}\

The hyperbolic Laplace operator $\gD:=-y^{2}(\dd_{xx}+\dd_{yy})$ acts (after unique extension) on the space
  $L^{2}(\G\bk \bH)$
  of
 square-integrable automorphic functions, and is self-adjoint and positive semi-definite.
Let $\gW=\gW(\G)\subset[0,\infty)$ denote the spectrum of $\gD$. The assumption that $\G$ has at least one cusp implies
the existence of continuous spectrum above $1/4$, that is, $[1/4,\infty)\subset\gW$ (there may also be embedded discrete spectrum in this range, which only occurs when $\G$ is a lattice  \cite{Selberg1956, Patterson1975}). Below $1/4$ the spectrum is finite \cite{LaxPhillips1982} and nonempty (by \eqref{eq:del12}); we denote these eigenvalues, often referred to as the ``exceptional spectrum,'' by
$$
0\le\gl_{0}<\gl_{1}\le\cdots\le\gl_{max}<\frac14,
$$
and introduce spectral parameters $1/2<s_{j}\le 1$ defined by
$$
\gl_{j} \  = \ s_{j}(1-s_{j}),
$$
so that
\be\label{eq:specSs}
1\ge s_{0}>s_{1}\ge \cdots\ge s_{max}>\frac12.
\ee
The bottom eigenvalue $\gl_{0}$ is simple, and is related to the geometry of $\gL$ via the Patterson-Sullivan formula \cite{Patterson1976, Sullivan1984}
$$
\gl_{0}\ = \ \gd (1-\gd),
$$
that is, $s_{0}=\gd$. 

\subsection{Algebra}\

We will use standard notation for the subgroups $N$, $A$, and $K$ of $G$, given by:
\be\label{eq:NAK}
N:=\mattwo1\R01,\quad A:=\left\{\diag( a,1/a):a>0\right\},\quad K:=\SO(2),
\ee
and containing typical elements
$$
n_{x}:=\mattwo 1x{}1, \quad
 a_{y}:=\mattwo{\sqrt y}{}{}{1/\sqrt y},
\quad k_{\gt}=\mattwo{\cos\gt}{\sin\gt}{-\sin\gt}{\cos\gt}.
$$
As right actions, they correspond, respectively, to the unipotent flow, geodesic flow, and
rotation of the tangent vector, keeping the base point fixed.
On the other hand, as left actions, the correspond, respectively, to horizontal translation, scaling, and
motion around a hyperbolic circle centered at $i$.
 Haar measure $dg$ in Iwasawa coordinates $g=n_{x}a_{y}k_{\gt}$ is then given by $dg=dx\,y^{-2}dy\,d\gt$.
The
right-action by the
 semigroup $A^{+}:=\{a_{y}:y>1\}$ corresponds to the positive geodesic flow, so that a given point $\bx_{0}\in G\cong T^{1}\bH$ gives rise to the geodesic ray $\bx_{0}A^{+}$.

\subsection{Representation Theory}\

By the Duality Theorem \cite{GelfandGraevPS1966},
the spectral decomposition \eqref{eq:specSs} corresponds to
the decomposition of the right regular representation of $G$ on $L^{2}(\G\bk G)$ as
$$
L^{2}(\G\bk G)  \  = \ V_{0}\oplus V_{1}\oplus\cdots\oplus V_{max}\oplus V_{temp}.
$$
Here $V_{temp}$ consists of the tempered spectrum (a reducible subspace);
each $V_{j}$, $j=1,\dots,max$ is 
an irreducible
complementary series representation of parameter $s_{j}$;
and
$V_{0}$ is either the trivial representation (if $\G$ is a lattice), or  a complementary series representation of parameter $s_{0}=\gd$ (if $\G$ is thin).

We record here a Sobolev-norm version  \cite{BernsteinReznikov1998} of the exponential decay of matrix coefficients. Fix a basis $\sB=\{X_{1},X_{2},X_{3}\}$ for the Lie algebra $\fg$ of $G$, and given a smooth test function $\Psi\in C^{\infty}(\G\bk G)$, define the ``$L^{p}$, order-$d$'' Sobolev norm $\cS_{p,d}(\Psi)$  as
$$
\cS_{p,d}(\Psi) \ : = \ \sum_{\ord(\sD)\le d}\|\sD\Psi\|_{L^{p}(\G\bk G)}
.
$$
Here $\sD$ ranges over monomials in $\sB$ of order at most $d$.
\begin{thm}[\cite{CowlingHaagerupHowe1988, Shalom2000, Venkatesh2010}]
Let $(\pi,V)$ be a unitary $G$-representation, and assume there is a number $\gT>1/2$ so that  $V$ does not weakly contain any complementary series with parameter $s>\gT$. Then for all smooth $v,w\in V^{\infty}$, we have
\be\label{eq:HoweMoore}
|\<\pi(g).v,w\>| \ \ll \
\|g\|^{2(1-\gT)}
\cS_{2,1}(v)\cS_{2,1}(w)
,
\ee
as $\|g\|^{2}:=\tr(g\,{}^{\top}\!g)\to\infty$. The implied constant is absolute.
\end{thm}

Later we will encounter other Sobolev norms which are convex combinations of those above.

\subsection{Uniform Spectral Gaps}\label{sec:unif}\

Recalling the spectral decomposition \eqref{eq:specSs},
we call a number $\gT\in(1/2,\gd)$ a {\bf spectral gap}  for $\G$ if $\gT>s_{1}$. To make sense of a {\it uniform} such gap,
we 
assume integrality.
As in \secref{sec:introCosets}, if $\G$ consists of integer matrices, $\G<\PSL_{2}(\Z)$,
we may, given an integer $q\ge1$, define the level-$q$ ``congruence'' subgroup
$$
\G(q) \ := \ \ker(\G\to\PSL_{2}(\Z/q\Z)).
$$
Let $\gW(q)$ be the spectrum of $\gD$ on $L^{2}(\G(q)\bk \bH)$; clearly $\gW\subset\gW(q)$, and in general this inclusion is strict. We will call a number $\gT\in(1/2,\gd)$  a {\bf uniform spectral gap} for $\G$ if, for all $q\ge1$,
\be\label{eq:specGap}
\gW(q)\cap(\gd(1-\gd),\gT(1-\gT)] \ = \
\O,
\ee
that is, there are no  eigenvalues at any level $q$ in a neighborhood of the base eigenvalue $\gl_{0}=\gd(1-\gd)$. (Note that this definition is different from other related definitions in the literature.)
In a number of statements below, several quantities depend on the {\it spectral gap} in the lattice case, but only on the {\it critical exponent} in the thin case; to unify the two notions, we will say that such quantities depend on the {\bf first non-zero eigenvalue} of $\G$.

\subsection{Eisenstein Series}\

In this section we recall some basic facts from the theory of Eisenstein series. We will assume here that the Eisenstein series is with respect to a cusp at $\infty$ and note that Eisenstein series corresponding to other cusps are defined similarly, after conjugating that cusp to $\infty$.
In our applications, we will not have the flexibility to demand the cusp have width $1$, so deal below with arbitrary width.

The Eisenstein series corresponding to a cusp at $\infty$  of width $\gw>0$ is defined in the half-plane $\Re(s)>1$ by the convergent series
\be \label{eq:EisDef}
E(z,s) \ : = \
\frac{1}{\gw}\!
\sum_{\g\in \G_\infty\bk\G} \Im(\g z)^s,
\ee
where $\G_\infty=\left(\begin{smallmatrix} 1 & \omega \Z\\ 0& 1\end{smallmatrix}\right)$. 

Assume first that $\G$ is a lattice. Then
$E(z,s)$ has a meromorphic continuation to $\C$ with a functional equation 
sending $s\mapsto1-s$.
In fact, it is analytic in the half plane $\Re(s)> 1/2$ except for a simple pole at $s=1$ and perhaps finitely many poles 
\be\label{eq:residuals}
1>\sigma_1\ge\cdots\ge\sigma_h>1/2
.
\ee
These poles comprise the ``residual spectrum,'' which is a subset of the ``exceptional spectrum'' in \eqref{eq:specSs}; the remaining spectrum in this range, if any, is cuspidal. 
The residue at $s=1$ is 
$$
\Res_{s=1}E(s,z)={
1
\over \vol(\G\bk \bH)}
$$ 
and 
$$
\vf_{\sigma_j}(z)=\Res_{s=\sigma_j}E(s,z)
$$ are the residual forms.

For any integer $n\in \Z$ we also define the weight $2n$ Eisenstein series by
\be \label{eq:EisnDef}
E_n(z,s)\ := \ 
\frac{1}{\gw}\!
\sum_{\g\in \G_\infty\bk\G} \Im(\g z)^s \gep_\g(z)^{2n},
\ee
where
$$
\gep_{g}(z)={cz+d\over|cz+d|},\qquad g=\mattwo abcd.
$$
Unless the weight $n=0$, the $E_{n}$'s are all regular at $s=1$, that is, $E_{0}=E$ is the only Eisenstein series with a pole at $s=1$. In the range $\tfrac12<\Re(s)<1$, the poles of $E_{n}$ are the same $s=\gs_{1},\dots,\gs_{h}$ as those of $E$.
%
For each such pole $\gs=\sigma_j$ we denote by
\be\label{eq:vfres}\vf_{\sigma,n}=\Res_{s=\sigma}E_n(z,s),\ee the (un-normalized) residual form of weight $2n$.

We note for future reference that the weight-$2n$ Eisenstein series and the weight-$2n$ residual forms all lie in the space $(\G,2n)$ of functions on $\bH$ transforming by
\be\label{eq:gepTrans}
f(\g z)=\gep_{\g}(z)^{2n}f(z).
\ee
%
Still assuming that $\G$ is a lattice, we also have the following bounds
coming from the spectral decomposition of $L^2(\G\bk G)$, see, e.g.  \cite{Iwaniec1997book}. For any square-integrable $f\in (\G,2n)$,
we have the bound
\be \label{eq:L2bound}\frac{1}{2\pi}\int_\R |\<f,E_n(\cdot, \tfrac{1}{2}+ir)\>|^2dr\leq \norm{f}_2^2.\ee

When $\G$ is thin, we will only use the fact that the defining  series in \eqref{eq:EisDef} and \eqref{eq:EisnDef} converge absolutely in the range $\Re(s)>\gd$.

\subsection{Decay of Fourier Coefficients}\

In this subsection we wish to record the basic fact that the (parabolic) Fourier coefficients of an automorphic function decay in the cusp, in a uniform sense. The method we use to establish this is completely standard, 
though
the requisite uniformity does not seem to be in the literature; hence we give sketches of proofs for the reader's convenience.
We again
assume that $\G$ has a cusp at $\infty$ of width $\gw>0$, that is, the isotropy group $\G_{\infty}$ of $\infty$ is generated by the translation $z\mapsto z+\gw$.

Then 
a smooth, square-integrable,  $\G$-automorphic function $\Psi\in L^{2}\cap C^{\infty}(T^{1}(\G\bk \bH))$ 
has a Fourier expansion:
\be\label{eq:FourExp}
\Psi(x+iy,\gz) \ = \
\sum_{m\in\Z}
a_{\Psi}(m;y,\gz)
\
e_{\gw}(mx)
,
\ee
where $e_{\gw}(x):=e^{2\pi i x/\gw}$, and the Fourier coefficients are given by
\be\label{eq:aPsiIs}
a_{\Psi}(m;y,\gz)
\ :=\
\frac1\gw
\int_{0}^{\gw}
\Psi(x+iy,\gz)
e_{\gw}(-mx)dx
.
\ee
The next proposition records the decay of such as $y\to0$ (in a uniform statement);
the subscripts $F$
below stand for ``Fourier coefficients.''

\begin{prop}\label{prop:decayFC}
There is a ``Sobolev'' norm $\cS_{F}(\Psi)$ and constants $0< C_{F}<\infty$ and $0<\ga_{F}<1$, so that, uniformly over all $y>0$ and $m\in\Z\setminus\{0\}$, we have
\be\label{eq:decayFC}
\left|
a_{\Psi}(m;y,\gz)
\right|
\ll
\cS_{F}(\Psi)\ |m|^{C_{F}}\ y^{\ga_{F}}
.
\ee
The constants $\ga_{F}$ and $C_{F}$ depend on the
first non-zero eigenvalue of $\G$.
\end{prop}

The Sobolev norm and exact values of the constants $C_{F}$ and $\ga_{F}$ are given below in \eqref{eq:cSFis}, \eqref{eq:CFis}, and \eqref{eq:gaFis}, respectively;
the last claim of the proposition is
 then 
 clear, 
 namely,  that the constants depend on the spectral gap when $\G$ is a lattice, and only on the critical exponent when $\G$ is thin.
As we are not striving for optimal exponents (recall \rmkref{rmk:optimize}), we have chosen to suppress their precise values so as not to clutter the presentation.
Recall also our convention that  implied constants may depend at most on $\G$, unless stated otherwise.

The Proposition is an easy consequence of another standard fact, namely the equidistribution of (pieces of) ``low-lying'' horocycles; the subscripts $H$
below stand for ``Horocycle pieces.''

\begin{prop}\label{prop:horoPiece}
There is a ``Sobolev'' norm $\cS_{H}(\Psi)$ and
constants $C_{H}<\infty$ and
 $0<\ga_{H}<1$, so that, uniformly over all $y>0$ and
open intervals
 $\cI\subset(0,\gw)$, we have
\be\label{eq:horoPiece}
\frac1{|\cI|}
\int_{\cI}
\Psi(x+iy,\gz)dx
\ = \
{1\over \vol(\G\bk G)}
\int_{\G\bk G}\Psi\ dg
+
O
\Bigg(
\cS_{H}(\Psi)\
|\cI|^{-C_{H}}\
y^{\ga_{H}}
\Bigg)
.
\ee
This statement holds whether $\G$ is a lattice or not, with the interpretation that the first  term on the right-hand-side of \eqref{eq:horoPiece} vanishes in the thin case.
The constants
 $C_{H}$ and
 $\ga_{H}$
 depend 
 on the first non-zero eigenvalue of $\G$.
\end{prop}

Again, the norm and constants
 are detailed in \eqref{eq:cSHis},  \eqref{eq:CHis},
and  \eqref{eq:gaHis}, which we have suppressed in the interest of exposition.
Much stronger versions of \eqref{eq:horoPiece} exist in the literature (at least in the lattice case, for which see, e.g., \cite{Strombergsson2004, SarnakUbis2015}), but for the reader's convenience, we provide a quick

\pf[Proof of \propref{prop:horoPiece} {[Sketch]}]\

We may assume that $y<1$ (for  the statement is obviously true otherwise
), and
we may moreover assume that $\Psi$ is right-$K$-invariant (or replace $\Psi$ by $\tilde\Psi:=\pi(k_{\gt}^{-1})\Psi$, where $2\gt$ is the angle of $\gz$
measured counterclockwise from the vertical). Then the left hand side of \eqref{eq:horoPiece} is
$$
\sM \ := \
\frac1{|\cI|}
\int_{x\in\cI}
\Psi(n_{x}a_{y})dx
.
$$
Let $\rho$ be a smooth, non-negative function on $\R$ with support in $[-1,1]$ and $\int_{\R}\rho=1$. For $\eta>0$ to be chosen later, concentrate $\rho$ to $\rho_{\eta}(x):=\eta^{-1}\rho(x/\eta)$.
Write $\overline N:={}^{\top}\!N$, $\overline n_{x}:={}^{\top}\!n_{x}$ for the opposite horocyclic group and element. Then the multiplication map
$$
N \times A\times \overline N \ \to \
G \ : \ (n,a,\overline n)\mapsto na\overline n
$$
is bijective on
an open neighborhood of the origin. Define $\xi:
G\to \R_{\ge0}$, supported in such a neighborhood, via:
$$
\xi(n_{x}a_{t}\overline n_{s}) \ :=
 \
 c_{\xi}\cdot
\left(\rho_{\eta}\star\frac1{|\cI|}\bo_{\cI}\right)\!(x)\cdot \rho_{\eta}(\log t)\cdot \rho_{\eta}(s),
$$
where $\star$ denotes convolution, and $c_{\xi}\asymp1$ is a constant (independent of $\eta$) chosen so that
 $\int_{G}\xi=1$.
 Automorphize $\xi$ to
$$
\Xi(g) \ := \ \sum_{\g\in\G}\xi(\g g),
$$
which is a function on $\G\bk G$ with $\int_{\G\bk G}\Xi =1$. Finally, consider the matrix coefficient:
$$
\sC \ := \ \<\pi(a_{y}).\Psi,\Xi\>
,
$$
which we evaluate in two ways.
Using the decay of matrix coefficients \eqref{eq:HoweMoore}, we see immediately that
$$
\sC \ =
{1\over\vol(\G\bk G)}\int_{\G\bk G}\Psi
+
O
\bigg(
y^{1-\gT}
\cS_{2,1}(\Psi)
\cS_{2,1}(\Xi)
\bigg)
,
$$
where $\gT=s_{1}+\vep$ is a spectral gap in the lattice case, and $\gT=\gd+\vep$ in the thin case (in which case the ``main'' term vanishes).
It is easy to estimate, crudely, that
$$
\cS_{2,1}(\Xi)\ll |\cI|^{-1} \eta^{-3}
.
$$

For a second evaluation of $\sC$, unfold the inner product to obtain
$$
\sC \ = \ \int_{NA\overline N}\Psi(n_{x}a_{t}\overline n_{s}a_{y}) \xi(n_{x}a_{t}\overline n_{s}) d\overline n_{s}da_{t}dn_{x}
$$
The ``wavefront lemma''  (in this case, trivial) states that $ \overline n_{s}a_{y} = a_{y}\overline n_{sy}$, and we estimate
$$
\Psi(n_{x}a_{y}a_{t}\overline n_{sy})
\ = \
\Psi(n_{x}a_{y})
+O\bigg(\eta\ \cS_{\infty,1}(\Psi)\bigg)
.
$$
Hence
$$
\sC \ = \ \int_{\R}\Psi(n_{x}a_{t})
\left(\rho_{\eta}\star\frac1{|\cI|}\bo_{\cI}\right)\!(x)\ dx
+O\bigg(\eta\ \cS_{\infty,1}(\Psi)\bigg)
\ = \
\sM
+O\bigg(\eta\ \cS_{\infty,1}(\Psi)\bigg)
.
$$
Combining the errors and choosing $\eta=y^{(1-\gT)/4}\cS_{2,1}(\Psi)^{1/4}\,\cS_{\infty,1}(\Psi)^{-1/4}|\cI|^{-1/4}$ gives
\eqref{eq:horoPiece}
 with
\be\label{eq:cSHis}
\cS_{H}(\Psi)\ : = \ \cS_{2,1}(\Psi)^{1/4} \ \cS_{\infty,1}(\Psi)^{3/4}
,
\ee
\be\label{eq:CHis}
C_{H}\ : = \ 1/4
,
\ee
and
\be\label{eq:gaHis}
\ga_{H}\ : = \ \frac{1-\gT}4
,
\ee
as claimed.
\epf

Equipped with \propref{prop:horoPiece}, we may now give a quick

\pf[Proof of \propref{prop:decayFC}]\ 

Again we may 
assume
that
 $\Psi$ is right-$K$-invariant and $y<1$.
Let $J\ge1$ be an
integer
  parameter to be chosen later, and write
$$
a_{\Psi}(m;y)
\ =\
\frac1\gw
\sum_{j=0}^{J-1}
\int_{\gw j/J}^{\gw(j+1)/J}
\Psi(x+iy)
e_{\gw}(-mx)dx.
$$
On each short interval, we estimate $e_{\gw}(-mx) = e(-m j/J)+O(|m|/J),$
whence
$$
a_{\Psi}(m;y)
\ =\
\sum_{j=0}^{J-1}
e(-mj/J)
\frac1\gw
\int_{\gw j/J}^{\gw(j+1)/J}
\Psi(x+iy)
dx
\ + \
O(\|\Psi\|_{\infty}|m|/J).
$$
Now on each little integral, we apply the equidistribution of pieces of ``low-lying'' horocycles  in the form \eqref{eq:horoPiece}, that is,
$$
\frac1\gw
\int_{\gw j/J}^{\gw(j+1)/J}
\Psi(x+iy)
dx
=
\frac1J
{1\over  \vol(\G\bk G)}
\int_{\G\bk G}\Psi\ dg
+
O
\Bigg(
\cS_{H}(\Psi)\
y^{\ga_{H}}
J^{C_{H}-1}
\Bigg)
.
$$
Inserting this expression into $a_{\Psi}(m;y)$ and using $m\neq0$, the roots of unity cancel 
out, 
leaving only error terms:
$$
|a_{\Psi}(m;y)|
\ \ll \
\cS_{H}(\Psi)\
y^{\ga_{H}}
J^{C_{H}}
+
\|\Psi\|_{\infty}|m|\ J^{-1}
.
$$
Setting
$$
J\asymp
\left(
\cS_{H}(\Psi)^{-1}\
y^{-\ga_{H}}\
\|\Psi\|_{\infty}\
|m|
\right)^{1/(C_{H}+1)}
,
$$
we arrive at \eqref{eq:decayFC} with 
\be\label{eq:cSFis}
\cS_{F}(\Psi) \ := \ \cS_{H}(\Psi)^{1/(C_{H}+1)}\ \|\Psi\|_{\infty}^{C_{H}/(C_{H}+1)}
,
\ee
\be\label{eq:CFis}
C_{F} \  : = \ C_{H}/(C_{H}+1)
,
\ee
and
\be\label{eq:gaFis}
\ga_{F}  \ := \ga_{H}/(C_{H}+1)
.
\ee
This completes the proof.
\epf
\begin{rmk}
It should be noted that actually Propositions \ref{prop:horoPiece} and \ref{prop:decayFC} are {\it equivalent}, in the sense that one can also use the uniform decay of Fourier coefficients to prove a version of \eqref{prop:horoPiece} (though with possibly worse exponents). 
\end{rmk}
\begin{rmk}
In the thin case, the proof of \propref{prop:decayFC} can be made {\it much} simpler. Namely, one can first trivially bound the $m^{th}$ coefficient by the constant one, $|a_{\Psi}(m;y,\gz)|\le |a_{\Psi}(0;y,\gz)|$, and then use \eqref{eq:horoPiece} with $\cI=(0,\gw)$ to estimate the constant coefficient.
 (%
 Note though that if $\G$ is a lattice, then of course
 $a_{\Psi}(0;y,\gz)$ need not decay!)
\end{rmk}

\newpage

\section{Equidistribution of Shears}\label{sec:mainPf}

Recall our running assumption that $\G<G=\PSL_{2}(\R)$ is a geometrically finite, Zariski dense, discrete group with at least one cusp, and hence critical exponent $\gd$ exceeding $1/2$.
As in \eqref{eq:muTdef}, we will study the limit as $|T|\to\infty$ of the measures
$$
\mu_{T}(\Psi) \ : = \
\int_{a\in A^{+}}\Psi(\bx_{0}\cdot a \cdot \fs_{T}) da.
$$


To study the equidistribution of such, we need an appropriate space of test functions; in particular, we will require smoothness and at least polynomial decay at the cusp. Toward this end, for any cusp $\fa$ of $\G$ and integer $m\ge1$, we introduce the space
$$
\sP_{\fa}^m(\G\bk G) \ \subset \ L^{2}\cap C^{\infty}(\G\bk G)
$$
of smooth, square-integrable, automorphic functions with the following added property.
%
We will state it
in the case $\fa=\infty$;
for a general cusp $\fa$, 
conjugate 
$\fa$ to $\infty$ in the standard way.

We require that, for each $\Psi\in\sP_\infty^m(\G\bk G)$, there are constants $1\le C_{\Psi}<\infty$ and $0<\ga_{\Psi}
$, such that 
\be\label{eq:polyDecay}
\left|
{\dd^j\over \dd\theta^j}\Psi(na_{y}k_\theta)
\right| \ < \ C_{\Psi}\, y^{-\ga_{\Psi}}
,
\ee
holds 
uniformly for all $j\leq m$, $y>C_{\Psi}$, and all $n\in N$, $k\in K$. That is, after a certain point high up in the specified cusp, we have completely uniform polynomial decay in $\Psi$ its first $m$ derivatives in $\fk=$ Lie$(K)$.
 Note that we make no
demands on decay properties (beyond square-integrability) in any other non-compact regions (cusps or possibly flares) of $\G\bk G$ besides the specified  cusp $\fa$. Also note that the space $\sP_{\fa}^{m}$ is non-empty, since, e.g., it contains the subspace of smooth, compactly supported functions, or better yet, cusp forms.

Our main theorem, from which \thmref{thm:main1} 
follows
 immediately, is the following.
\begin{thm}\label{thm:equi}
Let
$\bx_{0}\cdot A^{+}$
be
a cuspidal geodesic ray  ending in a cusp $\fa$ of $\G$,
and let 
$\Psi\in\sP_{\fa}^2(\G\bk G)$
be
  a test function 
(i.e., assume  \eqref{eq:polyDecay}
is satisfied
for all
$j\le2$). 
%
Then there is a finite-order ``Sobolev'' norm $\cS(\Psi)$
(which depends on the constants $C_{\Psi}$ and $\ga_{\Psi}$ in \eqref{eq:polyDecay}),
 and an $\eta>0$ 
depending only on the first non-zero eigenvalue of $\G$,
so that: if $\G$ is a lattice, 
$$
\mu_{T}(\Psi) \ = \ \log |T| \ \mu_{\G\bk G}(\Psi) + \mu_{\widetilde{Eis}}(\Psi) + O(
\cS(\Psi) T^{-\eta}),
$$
%
and if $\G$ is thin, then
$$
\mu_{T}(\Psi) \ = \ \mu_{{Eis}}(\Psi) + O(
\cS(\Psi)T^{-\eta}),
$$
as 
$|T|\to\infty$.
Here
$\mu_{\G\bk G}(\Psi):=\vol(\G\bk G)^{-1}\int_{\G\bk G}\Psi$
is
the Haar probability measure,
$\mu_{\widetilde{Eis}}$ 
is the
``regularized Eisenstein'' 
distribution
given
in \eqref{eq:muEis}, 
and
$\mu_{Eis}$ 
is the
distribution
given in \eqref{eq:muEisThin}. 
\end{thm}

As a first simplification, we can immediately apply an auxiliary conjugation to move $\bx_{0}$ to the origin $e\cong(i,\uparrow)$, whence the cusp $\fa$ moves to $\infty$. Unfortunately,
we have
thus
exhausted
our 
free parameters, 
and cannot control the {\it width} of the resulting cusp, which we denote $\gw$; that is, the isotropy group $\G_{\infty}$ is generated by the translation $z\mapsto z+\gw$.

As outlined in the introduction, the proof of \thmref{thm:equi} now proceeds in two stages, as encapsulated in the following two theorems.

\begin{thm}[Equidistribution in the ``strip'' $\fS=\G_{\infty}\bk G$]\label{thm:stage1}
For a test function $\Psi\in\sP^{2}_{\infty}(\G\bk G)$, define the measure:
\be\label{eq:muTfSdef}
\mu_{T,
\fS
}(\Psi) \ := \
\frac1{\vol(\G_{\infty}\bk N)}
\int_{n\in\G_{\infty}\bk N}
\int_{a\in A^{+}}
\Psi(
n\, a \, a_{1/T})da\, dn
\ = \
\frac1\gw
\int_{0}^{\gw}
\int_{1/T}^{\infty}
\Psi(
n_{x}\, a_{y}){dy\over y}\, dx.
\ee
(Recall that in this context, $da_{y}$ is $dy/y$, \emph{not} $dy/y^{2}$.)
Then
there is a ``Sobolev'' norm $\cS_{\fS}(\Psi)$ and a constant $\ga_{\fS}>0$, defined in \eqref{eq:cSfSis} and \eqref{eq:gafSis}, respectively,  so that
\be\label{eq:stage1}
\mu_{T}(\Psi)
\ =
\
\mu_{T,
\fS
}(\Psi)
\ + \
O\bigg(
\cS_{\fS}(\Psi)
T^{-\ga_{\fS}}
\bigg)
,
\ee
 as $|T|\to\infty$. 
 Here $\ga_{\fS}$ only depends on the first non-zero eigenvalue of $\G$.
\end{thm}

Note that 
\thmref{thm:stage1}
makes no
distinction between
 whether $\G$ is a lattice or thin. This 
dichotomy
 is only evident in the second stage:

\begin{thm}[Eisenstein 
distributions]\label{thm:stage2}
Let $\Psi\in\sP^{2}_{\infty}(\G\bk G)$ as above.

{\bf Lattice Case:} If $\G$ is a lattice in $G$, then there is a
 distribution 
  $\mu_{\widetilde{Eis}}$ defined in \eqref{eq:muEis}, and 
``residual''
distributions
  $\mu_{\sigma_j}$ corresponding to \eqref{eq:residuals}  and defined in \eqref{eq:muRes}, so that:
\begin{eqnarray*}
\mu_{T,\fS}(\Psi) & =&\mu(\Psi)\log(T)+\mu_{\widetilde{Eis}}(\Psi)\\
&&+\sum_{j=1}^h{T^{\sigma_j-1}\over \sigma_j-1} \mu_{\gs_j}(\Psi)+O(\cS_{2,1}(\Psi)T^{-1/2}).
\end{eqnarray*}
%
%
{\bf Thin Case:} If $\G$ is thin in $G$, then there is a 
distribution
$\mu_{{Eis}}$ defined in \eqref{eq:muEisThin} so that:
\be\label{eq:thm2Thin}
\mu_{T,\fS}(\Psi) \ = \
\mu_{{Eis}} (\Psi) \ + \ O\bigg( \cS_{H}(\Psi)T^{-\ga_{H}}\bigg).
\ee
Here $\cS_{H}$ and $\ga_{H}$ are as in \propref{prop:horoPiece}.
\end{thm}

It is clear that \thmref{thm:equi} follows immediately from Theorems \ref{thm:stage1} and \ref{thm:stage2}.

\pagebreak

\subsection{Stage 1:
Proof of \thmref{thm:stage1}
}\

We proceed with a series of elementary lemmata.
Beginning with the definition \eqref{eq:muTdef}, we express $\mu_{T}$ in terms of coordinates in $T^{1}(\G\bk \bH)$:
\be\label{eq:mu1}
\mu_{T}(\Psi) \ : = \
\int_{1}^{\infty}\Psi\bigg({y T\over \sqrt{T^{2}+1}}+ i {y\over \sqrt{T^{2}+1}},\uparrow\bigg) {dy\over y}
.
\ee
All of our manipulations below will not affect the direction of the tangent vector, so we drop the $\uparrow$.
(Alternatively, pretend $\Psi$ is right-$K$-invariant.)

\begin{lem}
With $C_{\Psi}$ and $\ga_{\Psi}$ from \eqref{eq:polyDecay}, we let
\be\label{eq:Udef}
U> C_{\Psi}T
\ee
be a parameter to be chosen later in \eqref{eq:Uis}.
Then
\be\label{eq:Lem1}
\mu_{T}(\Psi)
\ = \
\sM_{1}
\ + \
O\bigg(
\|\Psi\|_{\infty} T^{-2}
+
C_{\Psi}\left(\frac TU\right)^{\ga_{\Psi}}
\bigg)
,
\ee
where
\be\label{eq:sM1def}
\sM_{1} \ := \
\int_{1}^{U}
\Psi\left(y+i\frac yT\right){dy\over y}
.
\ee
\end{lem}
\pf
From \eqref{eq:mu1}, make a change of variables $y\mapsto y\sqrt{T^{2}+1}/T$, and simplify to
$$
\mu_{T}(\Psi) \ = \
\int_{T/\sqrt{T^{2}+1}}^{\infty}\Psi\bigg({y }+ i {y\over T}\bigg) {dy\over y}
\ = \
\int_{1}^{\infty}\Psi\bigg({y }+ i {y\over T}\bigg) {dy\over y}
\
+
\
O\bigg(
\|\Psi\|_{\infty} T^{-2}
\bigg)
.
$$
With $U$ as in \eqref{eq:Udef}, break the range of integration $[1,\infty)=[1,U]\cup(U,\infty)$. On the latter range, apply \eqref{eq:polyDecay},
whence
 \eqref{eq:Lem1} follows.
\epf

Now we invoke  the Fourier expansion \eqref{eq:FourExp}. Define
$$
\Psi^{\perp}(x+iy) \ := \ \sum_{m\in\Z\setminus\{0\}} a_{\Psi}(m;y)\ e_{\gw}(mx),
$$
so that
\be\label{eq:Decomp}
\Psi(x+iy) \ = \ a_{\Psi}(0;y) + \Psi^{\perp}(x+iy).
\ee
Inserting \eqref{eq:Decomp} into  \eqref{eq:sM1def} splits $\sM_{1}$ into a ``main term'' and ``error'':
$$
\sM_{1} \ = \ \sM_{2} + \sE_{1},
$$
where
\be\label{eq:sM2is}
\sM_{2}\ := \
\int_{1}^{U}
a_{\Psi}\left(0;\frac yT\right){dy\over y}
,
\ee
and
\be\label{eq:sE1is}
\sE_{1}
\ := \
\int_{1}^{U}
\Psi^{\perp}\left(y+i\frac yT\right){dy\over y}
.
\ee

We first analyze $\sM_{2}$.

\begin{lem}
Recalling the measure $\mu_{T,\fS}$ in \eqref{eq:muTfSdef}, we have
\be\label{eq:Lem2}
\sM_{2} \ = \
\mu_{T,\fS}(\Psi) \ + \
O\bigg(
C_{\Psi}\left(\frac TU\right)^{\ga_{\Psi}}
\bigg)
\ee
\end{lem}
\pf
Inserting \eqref{eq:aPsiIs} into \eqref{eq:sM2is} gives
$$
\sM_{2}
\ = \
\int_{1}^{U}
\frac1\gw
\int_{0}^{\gw}
\Psi\left(x+i\frac yT\right)
dx
{dy\over y}
\ = \
\frac1\gw
\int_{0}^{\gw}
\int_{1/T}^{U/T}
\Psi\left(x+i y\right)
{dy\over y}
dx
.
$$
Extending the $y$ integral from $U/T$ to $\infty$ and applying \eqref{eq:polyDecay} again gives the claimed main and error terms in \eqref{eq:Lem2}.
\epf

Returning to $\sE_{1}$ in \eqref{eq:sE1is},
our next goal is to incorporate the Fourier expansion, via the following
\begin{lem}
Let
$$
\sE_{2}
\ := \
\sum_{u=1}^{U}
\frac1{u^{2}}
\sum_{m\in\Z\setminus\{0\}}
{\left|
a_{\Psi}\!\left(
m; \frac uT
\right)
\right|
\over
|m|}
.
$$
Then
\be\label{eq:Lem3}
|\sE_{1}| \ \ll \
\sE_{2}
 \  + \
\cS_{\infty,1}(\Psi)\ {\log U\over T}
.
\ee
\end{lem}

\pf
We first straighten out the sheared integral in \eqref{eq:sE1is} by breaking it  into sums:
$$
\sE_{1}
\ = \
\sum_{u=1}^{U-1}
\int_{u}^{u+1}
\Psi^{\perp}\!\left(y+i\frac yT\right)
{dy\over y}
.
$$
On each interval, estimate
$$
\Psi^{\perp}\!\left(y+i\frac yT\right)
\ = \
\Psi^{\perp}\!\left(y+i\frac uT\right)
\ + \
O
\left(
\cS_{\infty,1}(\Psi) \ \frac1T
\right)
,
$$
 and Fourier expand
 $$
\Psi^{\perp}\!\left(y+i\frac uT\right)
\ = \
\sum_{m\neq0}
a_{\Psi}\!
\left(m;\frac uT\right)
e_{\gw}(my)
.
$$
Thus
$$
\sE_{1}
\ = \
\sum_{u=1}^{U-1}
\sum_{m\neq0}
a_{\Psi}\!
\left(m;\frac uT\right)
\left[
\int_{u}^{u+1}
e_{\gw}(my)
{dy\over y}
\right]
\ + \
O
\left(
\cS_{\infty,1}(\Psi) \ \frac{\log U}T
\right)
.
$$
Inserting absolute values and estimating the bracketed
term
 by partial integration
gives \eqref{eq:Lem3}, as claimed.
\epf

Our final task is to estimate $\sE_{2}$; we cannot directly use the decay of Fourier coefficients \eqref{eq:decayFC} in the full range of $m$, so introduce a parameter $M$ to be chosen later, and decompose
$$
\sE_{2} \ = \ \sE_{\ge} + \sE_{<},
$$
where for $\sq\in\{\ge, <\}$,
$$
\sE_{\sq} \ := \
\sum_{u=1}^{U}
\frac1{u^{2}}
\sum_{0\neq|m|\ \sq \ M}
{\left|
a_{\Psi}\!\left(
m; \frac uT
\right)
\right|
\over
|m|}
.
$$
We first estimate the large range trivially.
\begin{lem}
\be\label{eq:Lem4}
\sE_{\ge} \ \ll \
\|\Psi\|_{\infty} \ M^{-1/2}
.
\ee
\end{lem}
\pf
Cauchy-Schwarz and Parseval give:
\beann
\sE_{\ge}
& \ll &
\sum_{u=1}^{U}
\frac1{u^{2}}
\left(
\sum_{|m|\ge M}
{\left|
a_{\Psi}\!\left(
m; \frac uT
\right)
\right|^{2}
}
\right)^{1/2}
\left(
\sum_{|m|\ge M}
{1
\over
|m|^{2}}
\right)^{1/2}
\\
&\ll&
\sum_{u=1}^{U}
\frac1{u^{2}}
\left(
\frac1\gw
\int_{0}^{\gw}
\left|\Psi\left(x+i\frac uT\right)\right|^{2}
dx
\right)^{1/2}
M^{-1/2}
,
\eeann
which can be estimated by \eqref{eq:Lem4}, as claimed.
\epf

Finally, we estimate the range of small $m$ using decay of Fourier coefficients. Note that this is the {\it only} part of the argument
involving any spectral theory; nevertheless, thanks to  the uniformity of \propref{prop:decayFC}, we do not at this stage 
perceive  
any difference between the lattice and thin cases.

\begin{lem}
Recalling the Sobolev norm $\cS_{F}$ and constants $C_{F}$ and $\ga_{F}$ from \propref{prop:decayFC}, we have
\be\label{eq:Lem5}
\sE_{<} \ \ll \
\cS_{F}(\Psi) \ T^{-\ga_{F}} \  M^{C_{F}}
.
\ee
\end{lem}
\pf
Applying \eqref{eq:decayFC} gives
$$
\sE_{<} \ = \
\sum_{u=1}^{U}
\frac1{u^{2}}
\sum_{1\le |m|<M}
{\left|
a_{\Psi}\!\left(
m; \frac uT
\right)
\right|
\over
|m|}
\ \ll
\
\sum_{u=1}^{U}
\frac1{u^{2}}
\sum_{1\le |m|<M}
\cS_{F}(\Psi)
|m|^{C_{F}-1}
{\left|
 \frac uT
\right|^{\ga_{F}}
}
,
$$
which is bounded as claimed in \eqref{eq:Lem5}.
\epf

\pf[Proof of \thmref{thm:stage1}]
This is now a simple matter of combining the above lemmata. To balance \eqref{eq:Lem4} and  \eqref{eq:Lem5}, set
$$
M\ = \
\left(
\|\Psi\|_{\infty}\
\cS_{F}(\Psi)^{-1}\
T^{\ga_{F}}
\right)^{1/(C_{F}+1/2)}
,
$$
for a net error, crudely, of
\be\label{eq:sE2err}
\sE_{2} \ = \
\sE_{\ge}+\sE_{<}
\ 
= 
\
O\left(
\max(\cS_{F}(\Psi), \cS_{\infty,1}(\Psi))\cdot
T^{-\ga_{F}/(2C_{F}+1)}
\right)
.
\ee
To balance the error in  \eqref{eq:Lem1} and \eqref{eq:Lem2}  with that of   \eqref{eq:Lem3}, we take $U$ to be some power of $T$, say,
\be\label{eq:Uis}
U=T^{1+1/\ga_{\Psi}},
\ee
assuming that $T$ is large enough for \eqref{eq:Udef} to be satisfied.
Then the errors in
 \eqref{eq:Lem2} 
 and
 \eqref{eq:Lem1}
are 
 $O(C_{\Psi}/T)$, and
 the second error term in \eqref{eq:Lem3} is $O((1+\frac1{\ga_{\Psi}})\cS_{\infty,1}\log T/T)$,
%
which subsumes
 $O(\|\Psi\|_{\infty}T^{-2})$ in \eqref{eq:Lem1}.
On (again, crudely) setting
\be\label{eq:cSfSis}
\cS_{\fS}(\Psi) \ := \
\left(C_{\Psi}+\frac1{\ga_{\Psi}}\right)\cdot
\max(\cS_{F}(\Psi), \cS_{\infty,1}(\Psi))
,
\ee
and
\be\label{eq:gafSis}
\ga_{\fS} \ : = \ \ga_{F}/(2C_{F}+1),
\ee
as in \eqref{eq:sE2err},
one can verify
 directly
 that the net error is as claimed  in \eqref{eq:stage1}.

This completes the proof.
\epf

\


\subsection{Stage 2:
Proof of \thmref{thm:stage2}
}\


We first give the proof in the thin case, as it is significantly easier.

\subsubsection{Assume $\G$ is thin in $G$}\

Returning to \eqref{eq:muTfSdef}, write 
$\mu_{T,\fS}(\Psi)$ as:
$$
\mu_{T,
\fS
}(\Psi) \ = \
\frac1\gw
\left(
\int_{0}^{\infty}
-
\int_{0}^{1/T}
\right)
\int_{0}^{\gw}
\Psi(
z
,\uparrow)\
y\
dz
\ =: \
\sT_{1}-\sT_{2},
$$
say. Here we have set $dz:=dx\,dy/y^{2}$. We bound $\sT_{2}$ by
\be\label{eq:sT2bnd}
|\sT_{2}| \ \le \
\int_{0}^{1/T}
\left|
a_{\Psi}(0;y,\uparrow)
\right|\
y\
{dy\over y^{2}}
\ \ll
\
\int_{0}^{1/T}
\cS_{H}(\Psi)y^{\ga_{H}}
y\
{dy\over y^{2}}
\ = \
\cS_{H}(\Psi)\
T^{-\ga_{H}}
,
\ee
where we applied \eqref{eq:horoPiece} (with $\cI=(0,\gw)$).

Recalling that $\G_{\infty}=\mattwos1{\gw\Z}{}1 \ < \ \G$,
 we next deal with
\be\label{eq:sT1is}
\sT_{1}
\ := \
\frac1\gw
\int_{\G_{\infty}\bk \bH}
\Psi(
z,\uparrow)\
y\
dz
.
\ee
Note that the integral converges absolutely; for $y\to\infty$, this is due to  \eqref{eq:polyDecay}, while for $y\to0$, we can again use  \eqref{eq:horoPiece}.

For ease of exposition,
it is convenient at this point to
 first assume that $\Psi$
 is right-$K$-invariant, that is,
\be\label{eq:PsiKinv}
\Psi(
z,\gz) \  =\  \Psi(
z
).
\ee
Below
we 
detail
the
modifications needed  to handle the general case.

Recall from \eqref{eq:EisDef} that
$$
E(z,s)  \ := \
\frac1\gw
\sum_{\g\in\G_{\infty}\bk \G}
\Im(\g z)^{s}
$$
is
the  Eisenstein series at $\infty$ of a cusp of width $\gw$. Note that the defining sum converges absolutely and uniformly on compacta in the range $\Re(s)>\gd$, since $\G$ is assumed to be a thin subgroup of $G$. In particular, $E(z,s)$ is regular at $s=1$.

Then, letting $\sF$ be a fixed fundamental domain for $\G\bk\bH$, we can ``re-fold'' and write \eqref{eq:sT1is} as
$$
\sT_{1}
\ = \
\frac1\gw
\sum_{\g\in\G_{\infty}\bk\G}
\int_{\g\sF}
\Psi(
z)\
y\
dz
\ = \
\frac1\gw
\sum_{\g\in\G_{\infty}\bk\G}
\int_{\sF}
\Psi(
z)\
\Im(\g z)\
dz
\ = \
\<\Psi,E(\cdot,1)\>
.
$$
Setting
\be\label{eq:muEisThinK}
\mu_{Eis}(\Psi) \ := \ \<\Psi,E(\cdot,1)\>,
\ee
we immediately see that $\sT_{1}=\mu_{Eis}(\Psi)$, which combined with \eqref{eq:sT2bnd} gives:
$$
\mu_{T,
\fS
}(\Psi) \ = \
\mu_{Eis}(\Psi)
+O
\left(
\cS_{H}(\Psi)\
T^{-\ga_{H}}
\right)
,
$$
as claimed.

Finally, we remove the assumption \eqref{eq:PsiKinv} and 
extend
the proof 
to
the general case.
For a unit tangent vector $\gz$ at $z$, write $\gt\in[-\pi,\pi)$ for the ``angle''
of $\gz=\gz_{\gt}$, measured from the vertical $\uparrow$ counterclockwise.
We first decompose $\Psi(z,\gz)$ in a Fourier series in $\gz$, writing:
\be\label{eq:PsiFE}
\Psi(z,\gz) \ = \
\sum_{n\in\Z} \widehat\Psi_{n}(z) \ \chi_{n}(\gz),
\ee
where $\chi_{n}(\gz_{\gt})=e^{in\gt}$ in the above notation, and
$$
\widehat\Psi_{n}(z)
 \ := \
 \frac1{2\pi}
\int
_{-\pi}^{\pi}
\Psi(z,\gz_{\gt}) \
\overline\chi_{n}(\gz_{\gt})\
d\gt.
$$
Note that each $\hat\Psi_{n}$ lives in the space $(\G,2n)$ of functions on $\bH$ given in \eqref{eq:gepTrans}.

Returning to $\sT_{1}$ in \eqref{eq:sT1is}, we insert \eqref{eq:PsiFE} (with $\gz=\ \uparrow$), and 
``re-fold
'' 
again, obtaining:
\begin{eqnarray*}
\sT_{1}
& = &
\frac1\gw
\sum_{\g\in\G_{\infty}\bk\G}
\int_{\g\sF}
\sum_{n\in\Z}
\widehat\Psi_{n}(
z)\
\Im(z)\
dz\\
& = &
\frac1\gw
\sum_{n}
\sum_{\g\in\G_{\infty}\bk\G}
\int_{\sF}
\widehat\Psi_{n}(
z)\
\gep_{\g}(z)^{2n}
\Im(\g z)\
dz\\
& = &
\sum_n \<\widehat\Psi_n,E_n(\cdot,1)\>.
\end{eqnarray*}
Here, $E_{2n}(z,s)$ are the ``wight-$2n$'' Eisenstein series given by the series \eqref{eq:EisnDef};
these all converge absolutely for $\Re(s)>\gd$. 
The
absolute convergence of the sum over $n$ is guaranteed by  \eqref{eq:polyDecay} after taking two derivatives in $\gt$ and noting that ${{\dd^2\Psi}\over{\dd\gt^2}}\in \sP_{\fa}^0(\G\bk G)$. Then,
on defining 
\be\label{eq:muEisThin}
\mu_{Eis}(\Psi) \ := \ \sum_n \<\widehat\Psi_n,E_n(\cdot,1)\>,
\ee
 \eqref{eq:thm2Thin} follows immediately.
 Note that if $\Psi$ is $K$-invariant, the two definitions \eqref{eq:muEisThinK} and \eqref{eq:muEisThin} agree, and moreover $\mu_{Eis}$ is actually a measure. 
In general, $\mu_{Eis}$ is a
 distribution, as we need some derivatives of $\hat \Psi_{n}$ to ensure the convergence of \eqref{eq:muEisThin}.
This
 completes the proof in the thin case.

\ 

\subsubsection{Case $\G$ is a lattice in $G$}\

In this case, our analysis precedes in a similar fashion to that in \cite{Sarnak1981}.
We begin with the following
\begin{lem}
For $1<\gs<1+\alpha_\Psi$, we have
\be\label{eq:Lem12}
\mu_{T,\fS}(\Psi) \ = \ \sum_n
{1\over 2\pi i}\int_{(\gs)}{T^{s-1}\over s-1}
\<\widehat{\Psi}_n,E_n(\cdot,\bar s)\>
  ds
.
\ee
\end{lem}


\pf
Starting with \eqref{eq:muTfSdef}, write
\be\label{eq:muhT}
\mu_{T,\fS}(\Psi) \ = \
\frac1\gw
\int_{0}^{\infty}
\int_{0}^{\gw}
\Psi(z,\uparrow) \
h_{T}(y)\
dz
,
\ee
where we have set
$$
h_{T}(y) \ := \ y\cdot\bo_{\{y>1/T\}}.
$$
Note the Mellin transform/inverse pair:
$$
\widetilde{h_{T}}(s) \ := \ \int_{0}^{\infty}h_{T}(y)\ y^{-s}\  {dy\over y} \ = \ {T^{s-1}\over s-1},
$$
and
$$
h_{T}(y) \ = \
{1\over 2\pi i}\int_{(\gs)}{T^{s-1}\over s-1} \ y^{s} \ ds
.
$$
The first integral converges absolutely  for $\Re(s)=\gs>1$; the second 
is henceforth interpreted (after partial integration) as the absolutely convergent integral
\be\label{eq:hTy}
h_{T}(y) \ = \
{1\over 2\pi i}\int_{(\gs)}{T^{s-1}\over \log(Ty)(s-1)^{2}} \ y^{s} \ ds
.
\ee
Inserting \eqref{eq:hTy} into \eqref{eq:muhT} with the above convention gives
\be\label{eq:Lem11}
\mu_{T,\fS}(\Psi) \ = \
\frac1\gw
{1\over 2\pi i}\int_{(\gs)}{T^{s-1}\over s-1}
\int_{0}^{\infty}
\int_{0}^{\gw}
\Psi(z ,\uparrow
) \
 y^{s}
\ 
dz
  \ ds
,
\ee
which is absolutely convergent in the range $1<\gs<1+\ga_{\Psi}$ using
 \eqref{eq:polyDecay}. 

Now we proceed as in the thin case,
decomposing
$$\Psi(z,\uparrow)\  = \ \sum_n\widehat{\Psi}_n(z)$$ 
and
``unfolding''; for each $n\in \Z$ this gives
$$
 \int_{\G_{\infty}\bk\bH} \widehat{\Psi}_n(z)  \ y^{ s} \  dz 
\ =\ 
 \<\widehat{\Psi}_n,E_n(\cdot,\bar s)\> 
.
$$
Summing over 
$n$ 
and inserting into \eqref{eq:Lem11}
gives \eqref{eq:Lem12}, as claimed.
\epf

To finish the proof  of \thmref{thm:stage2}, we make the following definition:
 $$
 \widetilde{E}_n(z,s)=\left\{\begin{array}{ll} E_n(z,s) & n\neq 0\\ E(z,s)-\frac{1}{\vol(\G\bk \bH)(s-1)}& n=0\end{array}\right.,$$
 which, again, is regular at $s=1$ for all $n$.
Then \eqref{eq:Lem12} can be rewritten as
\bea
\mu_{T,\fS}(\Psi) 
&=&
\mu_{\G\bk G}(\Psi)\log(T)
\ + \ \sum_{n\in \Z}
{1\over 2\pi i}\int_{(\gs)}{T^{s-1}\over s-1}
\<\widehat{\Psi}_n,\widetilde{E}_n(\cdot,\bar s)\>ds
,
\eea
where we used that
$$
\<\widehat\Psi_{0},{1\over \vol(\G\bk \bH)} \>\ = \ \mu_{\G\bk G}(\Psi),
\qquad\text{and}\qquad
{1\over 2\pi i}\int_{(\gs)}{T^{s-1}\over (s-1)^{2}}
ds
\ = \
\log T.
$$
Now shifting the contour of integration to $\Re(s)=\tfrac{1}{2}$, we pick up
residues from the simple pole at $s=1$ and from the residual spectrum at $s=\gs_{j}$ as in \eqref{eq:residuals}.
The residue at $s=1$ is
\be \label{eq:muEis}
\sum_{n\in \Z}
\<\widehat{\Psi}_n,\widetilde{E}_n(\cdot,1)\> \ =:\ \mu_{\widetilde{Eis}}(\Psi),
\ee
that is, this is our ``second-order'' contribution, and is a distribution (as opposed to a measure) since $\Psi$ is not assumed to be $K$-finite.
Note that if $\Psi$ is $K$-fixed, then \eqref{eq:muEis} simplifies to just
%
\be \label{eq:muEisK}
\mu_{\widetilde{Eis}}(\Psi)
 \ = \
\<\Psi,\widetilde{E}_{0}(\cdot,1)\> ,
\ee
as claimed in \eqref{eq:tilEisIs}.

Each pole $s=\gs_{j}$ contributes
 a residue ${T^{\sigma_j-1}\over \sigma_j-1} \mu_{\sigma_j}(\Psi)$, where
\be \label{eq:muRes}\mu_{\sigma_j}(\Psi):= \sum_{n\in \Z}
\<\widehat{\Psi}_n,\vf_{\sigma_j,n}\>,
\ee
with $\vf_{\sigma_j,n}$ the ``weight-$2n$'' residual form given in \eqref{eq:vfres}.
Note that these distributions are exactly the same as those arising in Sarnak's analysis \cite[p. 737]{Sarnak1981}.

We thus obtain 
\begin{eqnarray*}
\mu_{T,\fS}(\Psi) & =&\mu(\Psi)\log(T)+\mu_{\widetilde{Eis}}(\Psi)
+\sum_{j=1}^h{T^{\sigma_j-1}\over \sigma_j-1} \mu_{\gs_j}(\Psi)\\
&&
\hskip.5in
+\sum_n {1\over 2\pi i}\int_{(1/2)}{T^{s-1}\over s-1}
\<\widehat{\Psi}_n,E_n(\cdot,\bar s)\> ds
.
\end{eqnarray*}
Finally, taking absolute values and
combining
Cauchy Schwartz
with \eqref{eq:L2bound},
 we can bound each of the terms in the last sum by
\begin{eqnarray*}
\left|
{1\over 2\pi i}\int_{(1/2)}{T^{s-1}\over s-1}
\<\widehat{\Psi}_n,E_n(\cdot,\bar s)\> ds
\right|
\ll T^{-1/2} \|\widehat{\Psi}_n\|_2.
\end{eqnarray*}
On estimating $\sum_n \|\widehat{\Psi}_n\|_2\ll \cS_{2,1}(\Psi)$,
we finally
 conclude the proof of \thmref{thm:stage2}.




\newpage

\section{Application 1: Moments of $L$-functions}\label{sec:thmMom}
Theorem \ref{thm:secMom} now follows readily from Theorem \ref{thm:main1}, as
we explain below.
Recall that we will illustrate the method on the simplest case of   $f$ being a weight-$k$ holomorphic Hecke cusp form on $\PSL_{2}(\Z)$; the calculation for general cuspidal $\GL(2)$ automorphic representations is similar.

Let
 $\Psi(x+iy)=|f(x+iy)|^{2}y^{k}$, 
 and use \eqref{eq:Hecke2} to
 write the left hand side of \eqref{eq:thmMom} as
\beann
\frac1{2\pi}
\int_{\R}
|L(f,\tfrac12+it)|^{2}|\cW_{k}(\tfrac12+it,T)|^{2}
dt
\ = \
\left(
\int_{0}^{
1/\widetilde T
}
+
\int_{
1/\widetilde T
}^{\infty}
\right)
\Psi\left(Ty+i{y}\right)
{dy\over y}
,
\eeann
where we have set 
$$
\widetilde T\ := \ \sqrt{T^{2}+1}
$$ 
for convenience.
\thmref{thm:main1} can be applied directly to the range $[1/\widetilde T,\infty)$, but the range $(0,1/\widetilde T)$ must be manipulated.
Changing variables $y\mapsto1/y$, using the automorphy of $\Psi$ that $\Psi(-1/z)=\Psi(z)$, and changing $y\mapsto \widetilde T^{2} y$ gives:
\beann
\int_{0}^{1/\widetilde  T}
\Psi\left(Ty+i{y}\right)
{dy\over y}
&=&
\int_{
\widetilde T}
^{\infty}
\Psi\left(\frac Ty+\frac i{y}\right)
{dy\over y}
\ = \
\int_{
\widetilde T}
^{\infty}
\Psi\left(-\frac {yT}{\widetilde T^{2}}
+
\frac {iy}{\widetilde T^{2}}
\right)
{dy\over y}
\\
&=&
\int_{
1/\widetilde T}
^{\infty}
\Psi\left(- {yT}
+
{iy}
\right)
{dy\over y}
.
\eeann
Now we can apply \thmref{thm:main1} to both contributions, giving
\be\label{eq:applyThm1}
\frac1{2\pi}
\int_{\R}
|L(f,\tfrac12+it)|^{2}|\cW_{k}(\tfrac12+it,T)|^{2}
dt
\ = \
2\mu_{\G\bk G}(\Psi)\log(T)+2\mu_{\widetilde{Eis}}(\Psi)+O_\Psi(T^{-\eta})
.
\ee
The first term is of course 
$$
\mu_{\G\bk G}(\Psi) \  =  \ \frac{\|f\|^2}{\vol(\G\bk\bH)}
,
$$ where the norm is with respect to the Petersson inner product. 
It remains to show that the second term, that is, the Eisenstein measure $\mu_{\widetilde{Eis}}(\Psi)$, can be expressed as 
special value (at the edge of the critical strip) of a symmetric square $L$-function. Note that $\Psi$ is a function on $\bH$, that is, as a function on $G$ it is right-$K$-invariant; therefore
$\mu_{\widetilde{Eis}}(\Psi)$ is determined by the simpler expression \eqref{eq:muEisK} (or \eqref{eq:tilEisIs}).

\begin{prop}\label{prop:tilE}
With the above notation, we have
$$
\mu_{\widetilde{Eis}}(\Psi) \ = \
{\|f 
\|^{2}\over \vol(\G\bk\bH)}
\left(
{\gL'\over \gL}(\sym^{2}f,1)
+
\g
-2
{\gz'\over \gz}(2)
\right)
.
$$
\end{prop}

Clearly \propref{prop:tilE} inserted into \eqref{eq:applyThm1} gives the right hand side of \eqref{eq:thmMom}, completing the proof of \thmref{thm:secMom}.

\pf[Proof of \propref{prop:tilE}]

To evaluate
$$
\cI \ := \
\mu_{\widetilde {Eis}}(
\Psi) \ 
= \
\<|f|^{2}y^{k}, \widetilde E(\cdot,1)\>
,
$$
 use \eqref{eq:tilEs} to write
\be\label{eq:cI1}
\cI \ = \
\lim_{s\to1}
\left(
\<
|f|^{2}y^{k}
,
\widetilde E(\cdot,\bar s)
\>
\right)
\ = \
\lim_{s\to1}
\left(
\<
|f|^{2}y^{k},
E(z,\bar s)
\>
-
{1\over (s-1)V}
\|f\|^{2}
\right)
,
\ee
where
$$
V \ := \
\vol(\G\bk \bH).
$$
We analyze $\<
|f|^{2}y^{k}
,E(z,s)\>
$ by standard Rankin-Selberg theory; more generally for  two cusp forms $f$ and $g$ of weight $k$, we have
\beann
\<f\,\bar g
\, y^{k},  E(\cdot,\bar s) \>
&=&
\int_{\G_{\infty}\bk\bH}f(z)\bar g(z) y^{k} y^{s}{dx}{dy\over y^{2}}
\\
&=&
\int_{0}^{\infty}
\int_{0}^{1}
\sum_{n\ge1}a_{f}(n)e^{2\pi i n x}e^{-2\pi ny}
\sum_{m\ge1}\overline{a_{g}(m)}e^{-2\pi i m x}e^{-2\pi my}
 y^{k} y^{s}{dx}{dy\over y^{2}}
\\
&=&
\sum_{n\ge1}
{a_{f}(n)
\overline{a_{g}(n)}
\over
n^{s+k-1}
}
\int_{0}^{\infty}
e^{-4\pi y}
 y^{k-1} y^{s}{dy\over y}
 \\
 &=&
(4\pi)^{-(s+k-1)}
\G(s+k-1)
 L(f\otimes \bar g,s) 
\  = \ \gL(f\otimes \bar g,s).
\eeann
When $f=g$, the Rankin-Selberg $L$-function factors (see, e.g., \cite[p. 232]{Iwaniec1997book}) 
as
$$
L(f\otimes \bar f,s) \ = \ {\gz(s)\over \gz(2s)} L (\sym^{2}f,s)
.
$$
Hence
\be\label{eq:RSsym2}
\<|f|^{2} y^{k}, E(\cdot,\bar s)\>
\ =\
 \gL(f\otimes\bar f,s)
\ = \
{\gz(s)\over \gz(2s)}
\gL(\sym^{2}f,s)
,
 \ee
 where
$ 
\gL(\sym^{2}f,s)
$ 
is
as in \eqref{eq:sym2f}.
Taking residues at $s=1$ on both sides of \eqref{eq:RSsym2} gives
\be\label{eq:f2ovV}
{\|f\|^{2} \over V}
\ = \
{1\over \gz(2)}
\gL(\sym^{2}f,1)
.
\ee
Inserting \eqref{eq:f2ovV} and \eqref{eq:RSsym2} into \eqref{eq:cI1} gives
$$
\cI \ = \
\lim_{s\to1}
\left(
{\gz(s)\over \gz(2s)}
\gL(\sym^{2}f,s)
-
{1\over (s-1)}
{1\over \gz(2)}
\gL(\sym^{2}f,1)
\right)
.
$$
Using $
\gz(s)-\frac1{s-1}\to \g
$
as $s\to1$ (Euler's constant), and elementary calculus,
 we have that
\beann
\cI&
=
&
\g{\gL(\sym^{2}f,1)\over \gz(2)}
+
{\gL'(\sym^{2}f,1)\over \gz(2)}
-
2{\gL(\sym^{2}f,1)\over \gz(2)}\cdot{\gz'\over\gz}(2)
\\
&
=&
{\|f\|^{2}\over V}
\left(
{\gL'\over \gL}(\sym^{2}f,1)
+
\g
-2
{\gz'\over \gz}(2)
\right)
,
\eeann
on using 
\eqref{eq:f2ovV} again.
This completes the proof.
\epf

Note that one could also extend our method to Eisenstein series, and then evaluate the (weighted) fourth moment of the Riemann zeta function using \thmref{thm:main1}. 

\subsection{Subconvexity?}\

We leave open the problem of extracting from the effective second moment \eqref{eq:thmMom} a subconvex bound 
$$
|L(f,\tfrac12+it)| \ \ll_{f}\ |t|^{1/2-\eta}
$$ 
in the $t$-aspect. Such is already known \cite{Good1986, PetridisSarnak2001} in the Maass case via trace formulae, explicit expansions, and shifted convolutions, but it would be interesting to give a new proof using only equidistribution. (Of course the general $\GL(2)$ subconvexity problem has been resolved \cite{MichelVenkatesh2010}; the interest here would be in the method used.)

The key issue is that the archimedean factor $|\cW_{k}|^{2}$ in \eqref{eq:thmMom}
is a smooth weight, which does not allow truncation; if the weight could be replaced by a sharp cutoff while still having a power savings rate, then the subconvexity bound would follow
immediately. This could be accomplished by finding a function $\Psi_{X}(T)$ so that
\be\label{eq:PsiX}
\int_{\R}\Psi_{X}(T)|\cW_{k}(\tfrac12+it,T)|^{2}dT \ \overset{?}{=} \ \bo_{|t|<X};
\ee
indeed, then one would multiply both sides of \eqref{eq:thmMom} by $\Psi_{X}(T)$ and integrate in $T$, 
obtaining
\beann
\frac1{2\pi}\int_{|t|<X}|L(f,\tfrac12+it)|^{2}dt
&=&
\int_{\R}\Psi_{X}(T)\left(C_{1}\log T + C_{2} +O(T^{-\eta})\right)dT \\ 
&\overset{?}{=}  &
C_{1}'X\log X+C_{2}'X + O(X^{1-\eta'}).
\eeann
Another approach is ``shorten the interval,'' that is, to replace the right hand side of \eqref{eq:PsiX} by $\bo_{|t-X|<Y}$, with $Y<X^{1-\eta}$.

Either way, one would need to invert the ``$\cW$-transform'': 
$$
\Psi(T)\ \mapsto \ \widetilde\Psi(t) \  := \ \int_{\R}\Psi(T)|\cW_{k}(\tfrac12+it,T)|^{2}dT.
$$
Unfortunately, there are basic difficulties with said inversion, namely a Paley-Weiner (or Heisenberg uncertainty) analysis shows that the transform has insufficient harmonics to be invertible
and functions $\Psi_{X}$ as above do not exist,
even in this simple holomorphic case! (Cf. the related discussion in, e.g., \cite[Appendix]{GhoshReznikovSarnak2013}.)
The case of non-holomorphic  Hecke-Maass forms is seemingly even more complicated as the weights \eqref{eq:cWk} will involve Bessel functions.

A potential method to circumvent this issue
(since our equidistribution theorem is proved in the generality of the unit tangent bundle)
 is to use all the harmonics afforded us by $f$, that is, by applying Maass raising and lowering operators.  This does not change the $L$-function, but results in effective second moments with 
a large span of
 weight functions $\cW$. One can hope that enough combinations of these can recover the desired sharp cutoff  functions $\Psi_{X}$, and we plan  to return to this question later.

\newpage

\section{Application 2: Counting and Non-Equidistribution}\label{sec:counting}

\subsection{
Proof of \thmref{thm:count}}\
 
As the method of counting from equidistribution is by now completely standard, we give a brief sketch (only the setup is not completely obvious).
Let $G=\SL_{2}(\R)$ be the spin double-cover $G\overset{\iota}\longrightarrow\SO_{Q}^{\circ}(\R)$ of the (identity component of the) special orthogonal group preserving an indefinite ternary quadratic form $Q$. Let $\G<G$ be discrete, Zariski-dense, geometrically finite, and have at least one cusp, and given $\bx_{0}\in\R^{3}$, let $\cO=\bx_{0}\,\iota(\G)$ be a discrete orbit. 
Let  $H=\Stab_{G}\bx_{0}$ be the stabilizer of $\bx_{0}$ in $G$, and let $\G_{H}:=\G\cap H$ be the stabilizer in $\G$.
Given an archimedean norm $\|\cdot\|$ on $\R^{3}$, we obtain a norm-$T$ ball $B_{T}$ in $H\bk G$ as in \eqref{eq:normiii}. Our goal is to estimate $\cN_{\cO}(T)=|\cO\cap B_{T}|$, that is,
$$
\cN_{\cO}(T) \ = \ \{\g\in\G_{H}\bk\G  \ : \  \|\bx_{0}\,\iota(\g)\|<T\}.
$$

Thanks to the discussion in \secref{sec:taxonomy} (see \tabref{tab:1}), there are only two new cases to prove, both occurring only when $H=\Stab_{G}\bx_{0}$ is diagonalizable.
We can choose the spin cover $\iota$ up to   conjugation, and hence can assume
 that $H=A$. 
Having made such a choice, we will henceforth drop $\iota$ from the notation.
To handle the two lacunary cases, we may assume that $\G_{H}$ is trivial. 
(In principle, $\G_{H}$ could be finite.)

We then break $\cN_{\cO}(T)$  into two contributions as follows. 
Recalling the shear $\fs_{t}$ in \eqref{eq:sTdef},
we  decompose each $g\in G=ANK$ uniquely as $g=a\fs_{t}k$, and write 
$$
G^{\pm}:=\{g=a\fs_{t}k\in G: a\in A^{\pm}\}
.
$$ 
Hence we can write 
$$
\cN_{\cO}(T)=\cN_{\cO}^{+}(T)+\cN_{\cO}^{-}(T)
,
$$ 
say, where
$$
\cN_{\cO}^{\pm}(T)  \ := \ \{\g\in\G\cap G^{\pm}  \ : \  \|\bx_{0}\,
\g
\|<T\}
,
$$
and treat only $\cN_{\cO}^{+}(T)$, the other contribution being the same (after conjugation).

If $\G$ is a lattice, then the ``lacunary''  case occurs only when both $0$ and $\infty$ (that is, the two endpoints of $A$) are cusps of $\G$. When $\G$ is thin,
the ``lacunary'' cases occur when
 at least one of $0$, $\infty$ is a cusp;
 \lemref{lem:cO} forces
  the other endpoint 
to be either a cusp or
in the free boundary. If $\infty$, say,  is in the free boundary, then $\cN_{\cO}^{+}(T)$ gives a contribution of order $N^{\gd}$, as described  below \eqref{eq:LaxP}. So
 to restrict attention to the lacunary case, we assume that $\infty$ is a cusp of $\G$. Now we continue with the
 standard smoothing/unsmoothing argument applied to the equidistribution theorem.
 For ease of exposition, assume that the norm $\|\cdot\|$ is right-$K$-invariant. (This assumption is standard to relax.)

Let $\psi:G\to\R_{\ge0}$ be a right-$K$-invariant bump function supported in an $\vep>0$ ball about the origin in $G/K$ with $\int_{G}\psi=1$. Set $\Psi(g):=\sum_{\g\in\G}\psi(\g g)$, so that $\int_{\G\bk G}\Psi=1$. Let
$$
f_{T}^{+}(g) := \bo_{\{\|\bx_{0} g\|<T, \ g\in G^{+}\}},
$$
and
$
\cF_{T}^{+}(g)\ := \ \sum_{\g\in\G}f_{T}^{+}(\g g ).
$
Then 
$$
\cF_{T}^{+}(e) \ = \ \cN_{\cO}^{+}(T)
$$
and
$$
\<\cF_{T}^{+},\Psi\> \ = \  \cN_{\cO}^{+}(T)(1+O(\vep))
,
$$
since $\Psi$ is a bump function about the origin.
Unfolding the inner product gives
$$
\<\cF_{T}^{+},\Psi\>
\ = \
\int_{G}f_{T}^{+}(g) \Psi(g) dg
\ = \
\int_{\fs_{t}}
f_{T}^{+}(\fs_{t})
\left[
\int_{a\in A^{+}}
 \Psi(a\fs_{t}) da
 \right]
 d\fs_{t}
 .
$$
Applying \thmref{thm:main1} to the bracketed term and integrating in $t$ completes the
 sketch
 of the two remaining cases of \thmref{thm:count}.
\\

\subsection{Proof of \propref{prop:noED}}\

As above, let
the stabilizer $H=A$ be diagonalizable,  let $\infty$ be a cusp of $\G$, and assume for ease of exposition that the norm $\|\cdot\|$ is right-$K$-invariant. 
The statement of \propref{prop:noED} assumes that
 $\G<\SL_{2}(\Z)$ is integral and thin.
For an integer $q\ge1$, let $\G(q)<\G$ be its level-$q$ principal congruence subgroup, and  
for a fixed $\vp\in\G/\G(q)$, let 
$$
\cO_{q,\vp} \ := \ \bx_{0}\vp\G(q)
$$ 
be  the congruence coset orbit. The corresponding counting function  is then
$$
\cN_{q,\vp}(T) \  := \  |\cO_{q,\vp}\cap B_{T}|.
$$ 
We claim that this count depends on $\vp$, that is, is not distributed uniformly among the cosets. 

One way to see
 this is to unravel the formalism of the previous proof, and note that $C_{1}$ in \eqref{eq:countThin} is essentially the evaluation at $s=1$ and some $z=z_{\vp}\in\bH$ of  the 
(unregularized) Eisenstein series $E(z,s)$ for $\G(q)$; there is no reason for these values to coincide for different $z_{\vp}$. An even easier way to see the non-equidistribution is to look at one picture.

\begin{figure}
\includegraphics[width=\textwidth]{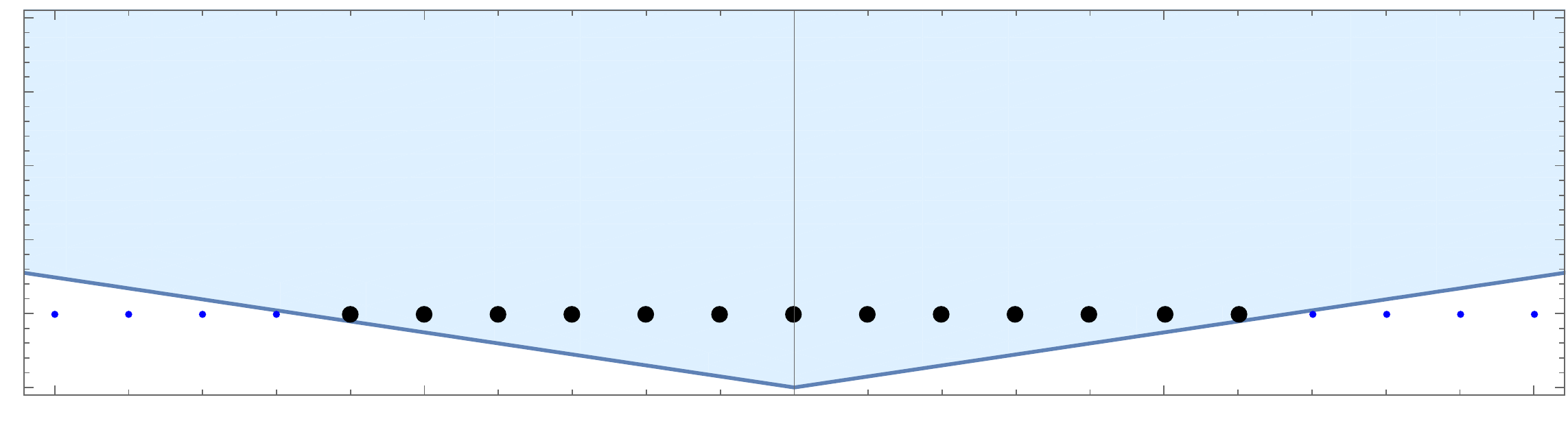}
\caption{The orbit $\bx_{0}\vp\G_{\infty,q}$ for $\vp=e$ inside $B_{T}\subset\bH$.}
\label{fig:nonED}
\end{figure}

Assume for simplicity that $\mattwos 1\Z{}1<\G$; then the isotropy group of $\infty$ in $\G(q)$ is 
$$
\G_{\infty,q} \ := \ \mattwos 1{q\Z}{}1 \ = \ \G(q)\cap \G_{\infty}.
$$ 
Certainly
the orbit $\cO_{q,\vp}$ 
 contains the points 
$$
\bx_{0}\vp\G_{\infty,q} \ \subset \ \cO_{q,\vp},
$$
so 
$$
\cN_{q,\vp}(T) 
\ \ge \  
|\bx_{0}\vp\G_{\infty,q}\cap B_{T}|
.
$$
Converting 
\figref{fig:BTA} from the disk $\bD$ to the hyperbolic plane $\bH$, we show in \figref{fig:nonED} how the shaded region $B_{T}$ contains the orbit points $\vp\G_{\infty,q}$ for $\vp=e$. 
A moment's 
reflection (or rather, translation)
shows that 
\eqref{eq:nonEDqT} holds for this orbit. 

Something similar would happen if one were to take congruence cosets with the subgroups $\G_{0}(q):=\{\mattwos abcd\in\G:c\equiv0(q)\}$, say, instead of $\G(q)$. The isotropy group $\G_{\infty}$ would now remain unchanged, but the same picture shows that, for the identity coset (with $\vp=e$), the number of points in an orbit is $\gg T$, whereas is average count should be of order $T/q$. 
We leave it as an interesting challenge to develop sieve methods which apply to this non-uniformly distributed (in the archimedean ordering) 
setting.

\newpage

\bibliographystyle{alpha}

\bibliography{AKbibliog}

\begin{thebibliography}{HKKL13}

\bibitem[BBKT14]{BlomerBrumleyKontorovichTemplier2014}
Valentin Blomer, Farrell Brumley, Alex Kontorovich, and Nicolas Templier.
\newblock Bounding hyperbolic and spherical coefficients of {M}aass forms.
\newblock {\em J. Th\'eor. Nombres Bordeaux}, 26(3):559--579, 2014.

\bibitem[Bea83]{Beardon1983}
Alan~F. Beardon.
\newblock {\em The Geometry of Discrete Groups}, volume~91 of {\em Graduate
  Texts in Mathematics}.
\newblock Springer-Verlag, New York, 1983.

\bibitem[BGS10]{BourgainGamburdSarnak2010}
Jean Bourgain, Alex Gamburd, and Peter Sarnak.
\newblock Affine linear sieve, expanders, and sum-product.
\newblock {\em Invent. Math.}, 179(3):559--644, 2010.

\bibitem[BGS11]{BourgainGamburdSarnak2011}
J.~Bourgain, A.~Gamburd, and P.~Sarnak.
\newblock Generalization of {S}elberg's 3/16th theorem and affine sieve.
\newblock {\em Acta Math}, 207:255--290, 2011.

\bibitem[BK14]{BourgainKontorovich2014b}
Jean Bourgain and Alex Kontorovich.
\newblock The affine sieve beyond expansion {I}: thin hypotenuses, 2014.
\newblock To appear, {\it IMRN}, {\tt arXiv:1307.3535}.

\bibitem[BKS10]{BourgainKontorovichSarnak2010}
J.~Bourgain, A.~Kontorovich, and P.~Sarnak.
\newblock Sector estimates for hyperbolic isometries.
\newblock {\em GAFA}, 20(5):1175--1200, 2010.

\bibitem[Blo08]{Blomer2008}
Valentin Blomer.
\newblock Sums of {H}ecke eigenvalues over values of quadratic polynomials.
\newblock {\em Int. Math. Res. Not. IMRN}, (16):Art. ID rnn059. 29, 2008.

\bibitem[BR98]{BernsteinReznikov1998}
Joseph Bernstein and Andr{\'e} Reznikov.
\newblock Sobolev norms of automorphic functionals and {F}ourier coefficients
  of cusp forms.
\newblock {\em C. R. Acad. Sci. Paris S\'er. I Math.}, 327(2):111--116, 1998.

\bibitem[CHH88]{CowlingHaagerupHowe1988}
M.~Cowling, U.~Haagerup, and R.~Howe.
\newblock Almost {$L\sp 2$} matrix coefficients.
\newblock {\em J. Reine Angew. Math.}, 387:97--110, 1988.

\bibitem[Del42]{Delsarte1942}
J.~Delsarte.
\newblock Sur le gitter fuchsien.
\newblock {\em C. R. Acad. Sci. Paris}, 214:147--179, 1942.

\bibitem[Del74]{Deligne1974}
Pierre Deligne.
\newblock La conjecture de {W}eil. {I}.
\newblock {\em Inst. Hautes \'Etudes Sci. Publ. Math.}, (43):273--307, 1974.

\bibitem[DG10]{DiaconuGarrett2010}
A.~Diaconu and P.~Garrett.
\newblock Subconvexity bounds for automorphic {$L$}-functions.
\newblock {\em J. Inst. Math. Jussieu}, 9(1):95--124, 2010.

\bibitem[DGG12]{DiaconuGarrettGoldfeld2012}
Adrian Diaconu, Paul Garrett, and Dorian Goldfeld.
\newblock Moments for {$L$}-functions for {$GL_r\times GL_{r-1}$}.
\newblock In {\em Contributions in analytic and algebraic number theory},
  volume~9 of {\em Springer Proc. Math.}, pages 197--227. Springer, New York,
  2012.

\bibitem[Dol98]{Dolgopyat1998}
Dmitry Dolgopyat.
\newblock On decay of correlations in {A}nosov flows.
\newblock {\em Ann. of Math. (2)}, 147(2):357--390, 1998.

\bibitem[DRS93]{DukeRudnickSarnak1993}
W.~Duke, Z.~Rudnick, and P.~Sarnak.
\newblock Density of integer points on affine homogeneous varieties.
\newblock {\em Duke Math. J.}, 71(1):143--179, 1993.

\bibitem[EM93]{EskinMcMullen1993}
A.~Eskin and C.~McMullen.
\newblock Mixing, counting and equidistribution in lie groups.
\newblock {\em Duke Math. J.}, 71:143--180, 1993.

\bibitem[GGPS66]{GelfandGraevPS1966}
I.~M. Gelfand, M.~I. Graev, and I.~I. Pjateckii-Shapiro.
\newblock {\em Teoriya predstavlenii i avtomorfnye funktsii}.
\newblock Generalized functions, No. 6. Izdat. ``Nauka'', Moscow, 1966.

\bibitem[GJ78]{GelbartJacquet1978}
Stephen Gelbart and Herv{\'e} Jacquet.
\newblock A relation between automorphic representations of {${\rm GL}(2)$} and
  {${\rm GL}(3)$}.
\newblock {\em Ann. Sci. \'Ecole Norm. Sup. (4)}, 11(4):471--542, 1978.

\bibitem[Goo86]{Good1986}
A.~Good.
\newblock The convolution method for {D}irichlet series.
\newblock In {\em The {S}elberg trace formula and related topics ({B}runswick,
  {M}aine, 1984)}, volume~53 of {\em Contemp. Math.}, pages 207--214. Amer.
  Math. Soc., Providence, RI, 1986.

\bibitem[GRS13]{GhoshReznikovSarnak2013}
Amit Ghosh, Andre Reznikov, and Peter Sarnak.
\newblock Nodal domains of {M}aass forms {I}.
\newblock {\em Geom. Funct. Anal.}, 23(5):1515--1568, 2013.

\bibitem[HK14]{HongKontorovich2014}
J.~Hong and A.~Kontorovich.
\newblock Almost prime coordinates for anisotropic and thin {P}ythagorean
  orbits, 2014.
\newblock To appear, {\it Israel J. Math.} {\tt arXiv:1401.4701}.

\bibitem[HKKL13]{HulseKiralKuanLim2013}
Thomas~A. Hulse, E.~Mehmet Kiral, Chan~Ieong Kuan, and Li-Mei Lim.
\newblock Counting square discriminants, 2013.
\newblock Preprint, {\tt arXiv:1307.6606}.

\bibitem[Hoo63]{Hooley1963}
Christopher Hooley.
\newblock On the number of divisors of a quadratic polynomial.
\newblock {\em Acta Math.}, 110:97--114, 1963.

\bibitem[Hub56]{Huber1956}
H.~Huber.
\newblock \"uber eine neue {K}lasse automorpher {F}unktionen und ein
  {G}itterpunktproblem in der hyperbolischen {E}bene. i.
\newblock {\em Comment. Math. Helv.}, 30:20--62, 1956.

\bibitem[IK04]{IwaniecKowalski}
Henryk Iwaniec and Emmanuel Kowalski.
\newblock {\em Analytic number theory}, volume~53 of {\em American Mathematical
  Society Colloquium Publications}.
\newblock American Mathematical Society, Providence, RI, 2004.

\bibitem[Iwa97]{Iwaniec1997book}
Henryk Iwaniec.
\newblock {\em Topics in classical automorphic forms}.
\newblock American Mathematical Society, Providence, RI, 1997.

\bibitem[KO12]{KontorovichOh2012}
A.~Kontorovich and H.~Oh.
\newblock Almost prime {P}ythagorean triples in thin orbits.
\newblock {\em J. reine angew. Math.}, 667:89--131, 2012.
\newblock {\tt arXiv:1001.0370}.

\bibitem[Kon09]{Kontorovich2009}
A.~Kontorovich.
\newblock The hyperbolic lattice point count in infinite volume with
  applications to sieves.
\newblock {\em Duke J. Math.}, 149(1):1--36, 2009.
\newblock {\tt arXiv:0712.1391}.

\bibitem[Kon14]{Kontorovich2014}
Alex Kontorovich.
\newblock Levels of distribution and the affine sieve.
\newblock {\em Ann. Fac. Sci. Toulouse Math. (6)}, 23(5):933--966, 2014.

\bibitem[KS03]{KimSarnak2003}
H.~Kim and P.~Sarnak.
\newblock Refined estimates towards the {R}amanujan and {S}elberg conjectures.
\newblock {\em J. Amer. Math. Soc.}, 16(1):175--181, 2003.

\bibitem[LP82]{LaxPhillips1982}
P.D. Lax and R.S. Phillips.
\newblock The asymptotic distribution of lattice points in {E}uclidean and
  non-{E}uclidean space.
\newblock {\em Journal of Functional Analysis}, 46:280--350, 1982.

\bibitem[LS10]{LiuSarnak2010}
Jianya Liu and Peter Sarnak.
\newblock Integral points on quadrics in three variables whose coordinates have
  few prime factors.
\newblock {\em Israel J. Math}, 178:393--426, 2010.

\bibitem[Mar04]{Margulis2004}
Grigoriy~A. Margulis.
\newblock {\em On some aspects of the theory of {A}nosov systems}.
\newblock Springer Monographs in Mathematics. Springer-Verlag, Berlin, 2004.
\newblock With a survey by Richard Sharp: Periodic orbits of hyperbolic flows,
  Translated from the Russian by Valentina Vladimirovna Szulikowska.

\bibitem[MO13]{MohammadiOh2013}
Amir Mohammadi and Hee Oh.
\newblock Matrix coefficients, counting and primes for orbits of geometrically
  finite groups, 2013.
\newblock To appear, {\it JEMS}.

\bibitem[MV10]{MichelVenkatesh2010}
Philippe Michel and Akshay Venkatesh.
\newblock The subconvexity problem for {${\rm GL}_2$}.
\newblock {\em Publ. Math. Inst. Hautes \'Etudes Sci.}, (111):171--271, 2010.

\bibitem[Nau05]{Naud2005}
Fr{\'e}d{\'e}ric Naud.
\newblock Expanding maps on {C}antor sets and analytic continuation of zeta
  functions.
\newblock {\em Ann. Sci. \'Ecole Norm. Sup. (4)}, 38(1):116--153, 2005.

\bibitem[NS10]{NevoSarnak2010}
Amos Nevo and Peter Sarnak.
\newblock Prime and almost prime integral points on principal homogeneous
  spaces.
\newblock {\em Acta Math.}, 205(2):361--402, 2010.

\bibitem[OS13]{OhShah2013}
Hee Oh and Nimish~A. Shah.
\newblock Equidistribution and counting for orbits of geometrically finite
  hyperbolic groups.
\newblock {\em J. Amer. Math. Soc.}, 26(2):511--562, 2013.

\bibitem[OS14]{OhShah2014}
Hee Oh and Nimish~A. Shah.
\newblock Limits of translates of divergent geodesics and integral points on
  one-sheeted hyperboloids.
\newblock {\em Israel J. Math.}, 199(2):915--931, 2014.

\bibitem[Pat75]{Patterson1975}
S.~J. Patterson.
\newblock The {L}aplacian operator on a {R}iemann surface.
\newblock {\em Compositio Math.}, 31(1):83--107, 1975.

\bibitem[Pat76]{Patterson1976}
S.J. Patterson.
\newblock The limit set of a {F}uchsian group.
\newblock {\em Acta Mathematica}, 136:241--273, 1976.

\bibitem[PS01]{PetridisSarnak2001}
Yiannis~N. Petridis and Peter Sarnak.
\newblock Quantum unique ergodicity for {${\rm SL}_2(\mathscr O)\backslash\bold
  H^3$} and estimates for {$L$}-functions.
\newblock {\em J. Evol. Equ.}, 1(3):277--290, 2001.
\newblock Dedicated to Ralph S. Phillips.

\bibitem[Sar81]{Sarnak1981}
Peter Sarnak.
\newblock Asymptotic behavior of periodic orbits of the horocycle flow and
  {E}isenstein series.
\newblock {\em Comm. Pure Appl. Math.}, 34(6):719--739, 1981.

\bibitem[Sar84]{Sarnak1984}
Peter Sarnak.
\newblock Additive number theory and {M}aass forms.
\newblock In {\em Number theory ({N}ew {Y}ork, 1982)}, volume 1052 of {\em
  Lecture Notes in Math.}, pages 286--309. Springer, Berlin, 1984.

\bibitem[Sar85]{Sarnak1985a}
Peter~C. Sarnak.
\newblock Fourth moments of grossencharakteren zeta functions.
\newblock {\em Comm. Pure and Applied Math.}, 38:167--178, 1985.

\bibitem[Sel56]{Selberg1956}
A.~Selberg.
\newblock Harmonic analysis and discontinuous groups in weakly symmetric
  {R}iemannian spaces with applications to {D}irichlet series.
\newblock {\em J. Indian Math. Soc. (N.S.)}, 20:47--87, 1956.

\bibitem[Sha00]{Shalom2000}
Yehuda Shalom.
\newblock Rigidity, unitary representations of semisimple groups, and
  fundamental groups of manifolds with rank one transformation group.
\newblock {\em Ann. of Math. (2)}, 152(1):113--182, 2000.

\bibitem[Str04]{Strombergsson2004}
A.~Strombergsson.
\newblock On the uniform equidistribution of long closed horocycles.
\newblock {\em Duke Math. J.}, 123:507--547, 2004.

\bibitem[SU15]{SarnakUbis2015}
Peter Sarnak and Adri{\'a}n Ubis.
\newblock The horocycle flow at prime times.
\newblock {\em J. Math. Pures Appl. (9)}, 103(2):575--618, 2015.

\bibitem[Sul84]{Sullivan1984}
D.~Sullivan.
\newblock Entropy, {H}ausdorff measures old and new, and limit sets of
  geometrically finite {K}leinian groups.
\newblock {\em Acta Math.}, 153(3-4):259--277, 1984.

\bibitem[Tit51]{Titchmarsh1951}
E.~C. Titchmarsh.
\newblock {\em The Theory of the {R}iemann Zeta-Function}.
\newblock Oxford Univ. Press, London/New York, Oxford, 1951.

\bibitem[TT13]{TemplierTsimerman2013}
Nicolas Templier and Jacob Tsimerman.
\newblock Non-split sums of coefficients of {$GL(2)$}-automorphic forms.
\newblock {\em Israel J. Math.}, 195(2):677--723, 2013.

\bibitem[Ven10]{Venkatesh2010}
Akshay Venkatesh.
\newblock Sparse equidistribution problems, period bounds and subconvexity.
\newblock {\em Ann. of Math. (2)}, 172(2):989--1094, 2010.

\bibitem[Zag81]{Zagier1981}
D.~Zagier.
\newblock Eisenstein series and the {R}iemann zeta function.
\newblock In {\em Automorphic forms, representation theory and arithmetic
  ({B}ombay, 1979)}, volume~10 of {\em Tata Inst. Fund. Res. Studies in Math.},
  pages 275--301. Tata Inst. Fundamental Res., Bombay, 1981.

\end{thebibliography}

\end{document}